\documentclass[11pt]{article}
\usepackage{amsfonts,calc,amsfonts,amsmath,graphicx,float,amsthm}
\usepackage[sort&compress,numbers]{natbib}
\usepackage[lmargin=30mm,rmargin=30mm,bmargin=30mm,tmargin=30mm]{geometry}
\newtheorem{observation}{Observation}
\newtheorem{proposition}{Proposition}
\newtheorem{theorem}{Theorem}
\newtheorem{lemma}{Lemma}
\newtheorem{corollary}{Corollary}
\newtheorem{claim}{Claim}

\newcommand{\thmlabel}[1]{\label{thm:#1}}
\newcommand{\thmref}[1]{Theorem~\ref{thm:#1}}
\newcommand{\twothmref}[2]{Theorems~\ref{thm:#1} and \ref{thm:#2}}

\newcommand{\lemlabel}[1]{\label{lem:#1}}
\newcommand{\lemref}[1]{Lemma~\ref{lem:#1}}

\newcommand{\eqnlabel}[1]{\label{eqn:#1}}

\newcommand{\Eqnref}[1]{Equation~\eqref{eqn:#1}}

\newcommand{\figlabel}[1]{\label{fig:#1}}
\newcommand{\figref}[1]{Figure~\ref{fig:#1}}

\newcommand{\seclabel}[1]{\label{sec:#1}}
\newcommand{\secref}[1]{Section~\ref{sec:#1}}

\newcommand{\manysecref}[2]{Sections~\ref{sec:#1}--\ref{sec:#2}}

\newcommand{\corlabel}[1]{\label{cor:#1}}

\newcommand{\proplabel}[1]{\label{prop:#1}}
\newcommand{\propref}[1]{Proposition~\ref{prop:#1}}
\newcommand{\twopropref}[2]{Propositions~\ref{prop:#1} and \ref{prop:#2}}

\newcommand{\tablabel}[1]{\label{tab:#1}}
\newcommand{\tabref}[1]{Table~\ref{tab:#1}}

\newcommand{\obslabel}[1]{\label{obs:#1}}
\newcommand{\obsref}[1]{Observation~\ref{obs:#1}}

\newcommand{\NP}{\ensuremath{\mathcal{NP}}}
\newcommand{\Oh}[1]{\ensuremath{\protect\mathcal{O}(#1)}}

\newcommand{\half}{\ensuremath{\protect\tfrac{1}{2}}}

\newcommand{\ceil}[1]{\ensuremath{\protect\lceil#1\rceil}}
\newcommand{\CEIL}[1]{\ensuremath{\protect\left\lceil#1\right\rceil}}

\newcommand{\floor}[1]{\ensuremath{\protect\lfloor#1\rfloor}}

\newcommand{\N}{\ensuremath{\mathbb{N}}}

\newcommand{\R}{\ensuremath{\mathbb{R}}}

\newcommand{\rr}{\ensuremath{r}}
\newcommand{\ii}{\ensuremath{i}}
\newcommand{\jj}{\ensuremath{j}}

\newcommand{\paran}[1]{\textup{(}#1\textup{)}}

\newcommand{\segc}[2]{\ensuremath{
\raisebox{1ex}{\tiny$\bullet$\hspace*{-0.1em}}
\overline{\hspace*{0.05em}#1#2\hspace*{0.05em}}
\raisebox{1ex}{\hspace*{-0.1em}\tiny$\bullet$}}}

\newcommand{\seg}[2]{\ensuremath{\overline{#1#2}}}
\newcommand{\lin}[2]{\ensuremath{\overleftrightarrow{#1#2}}}
\newcommand{\ray}[2]{\ensuremath{\overrightarrow{#1#2}}}
\newcommand{\oray}[2]{\ensuremath{\overleftarrow{#1#2}}}
\newcommand{\WEDGE}[2]{\ensuremath{\bigtriangledown(#1,#2)}}

\newcommand{\aaa}{\hspace*{1ex}\textup{(a)}}
\newcommand{\bbb}{\hspace*{1ex}\textup{(b)}}
\newcommand{\ccc}{\hspace*{1ex}\textup{(c)}}
\newcommand{\ddd}{\hspace*{1ex}\textup{(d)}}

\newcommand{\qi}{\hspace*{1ex}\textup{(i)}}
\newcommand{\qii}{\hspace*{1ex}\textup{(ii)}}

\newcommand{\ONE}{\hspace*{1ex}\textup{(1)}}
\newcommand{\TWO}{\hspace*{1ex}\textup{(2)}}
\newcommand{\THREE}{\hspace*{1ex}\textup{(3)}}

\newcommand{\POINT}{\hspace*{2ex}$\bullet$\hspace*{1ex}}

\newcommand{\G}{\ensuremath{\mathcal{G}}}
\newcommand{\X}{\ensuremath{\mathcal{S}}}
\newcommand{\M}{\ensuremath{\mathcal{M}}}
\newcommand{\D}{\ensuremath{\mathcal{D}}}
\newcommand{\T}{\ensuremath{\mathcal{T}}}
\newcommand{\E}{\ensuremath{\mathcal{E}}}

\newcommand{\bt}[2][]{\ensuremath{\textup{\textsf{bt}}_{#1}(#2)}}
\newcommand{\ba}[2][]{\ensuremath{\textup{\textsf{ba}}_{#1}(#2)}}
\newcommand{\bsa}[2][]{\ensuremath{\textup{\textsf{bsa}}_{#1}(#2)}}
\renewcommand{\tt}[2][]{\ensuremath{\theta_{#1}(#2)}}
\newcommand{\ot}[2][]{\ensuremath{\theta_{\textsf{o}}(#2)}}
\newcommand{\got}[2][]{\ensuremath{\overline{\theta_{\textsf{o}}}(#2)}}
\newcommand{\gt}[2][]{\ensuremath{\overline{\theta}_{#1}(#2)}}
\newcommand{\sa}[1]{\ensuremath{\textup{\textsf{sa}}(#1)}}
\newcommand{\arb}[1]{\ensuremath{\textup{\textsf{a}}(#1)}}
\newcommand{\ga}[1]{\ensuremath{\overline{\textup{\textsf{a}}}(#1)}}
\newcommand{\gsa}[1]{\ensuremath{\overline{\textup{\textsf{sa}}}(#1)}}
\newcommand{\as}[1]{\ensuremath{\protect{V_{\widehat{#1}}}}}

\renewcommand{\baselinestretch}{1.07}
\newcommand{\Figure}[4][htb]{
\begin{figure}[#1]
	\begin{center}#3\end{center}
	\caption{\figlabel{#2}#4}
\end{figure}}

\begin{document}

\title{\bf Graph Treewidth and\\ Geometric Thickness Parameters\footnotemark[1]}

\footnotetext[1]{This research was initiated in 2005 while both authors were in the School of Computer Science at McGill University, Montr\'eal, Canada. A preliminary version of this paper was published in the \emph{Proceedings of the 13th International Symposium on Graph Drawing} (GD '05), Lecture Notes in Computer Science 3843:129--140, Springer, 2006. The full version was  published in \emph{Discrete and Computational Geometry} 37.4:641--670, 2007. That version contained a false conjecture, which is corrected on page~\pageref{CHANGE} of this version. No other changes have been made. }

\author{Vida Dujmovi\'c\footnotemark[2] \qquad David R. Wood\footnotemark[3]}

\footnotetext[2]{Department of Mathematics and Statistics, McGill University, Montr{\'e}al, Canada. Partially supported by NSERC, Centre de Recherches Math\'ematiques (CRM), and Institut des Sciences Math\'ematiques (ISM). Now at: School of Computer Science and Electrical Engineering, University of Ottawa, Ottawa, Canada (\texttt{vida.dujmovic@uottawa.ca}). Research supported by NSERC. } 

\footnotetext[3]{Departament de Matem{\`a}tica Aplicada II, Universitat Polit{\`e}cnica de Catalunya, Barcelona, Spain. Supported by a Marie Curie Fellowship of the European Community under contract 023865, and by the projects MCYT-FEDER BFM2003-00368 and Gen.\ Cat 2001SGR00224. Now at: 
School of Mathematical Sciences, Monash University, Melbourne, Australia (\texttt{david.wood@monash.edu}). Research supported by the Australian Research Council.}

\date{March 30, 2005; Revised: \today}

\maketitle

\begin{abstract}
Consider a drawing of a graph $G$ in the plane such that crossing edges are coloured differently. The minimum number of colours, taken over all drawings of $G$, is the classical graph parameter \emph{thickness}. By restricting the edges to be straight, we obtain the \emph{geometric thickness}. By additionally restricting the vertices to be in convex position, we obtain the \emph{book thickness}. This paper studies the relationship between these parameters and treewidth. 

Our first main result states that for graphs of treewidth $k$, the maximum thickness and the maximum geometric thickness both equal $\ceil{k/2}$. This says that the lower bound for thickness can be matched by an upper bound, even in the more restrictive geometric setting. Our second main result states that for graphs of treewidth $k$, the maximum book thickness equals $k$ if $k\leq2$ and equals $k+1$ if $k\geq3$. This refutes a conjecture of Ganley and Heath [\emph{Discrete Appl. Math.}\ 109(3):215--221, 2001]. Analogous results are proved for outerthickness, arboricity, and star-arboricity.
\end{abstract}


\newpage
\tableofcontents

\newpage
\section{Introduction}

Partitions of the edge set of a graph into a small number of `nice' subgraphs are in the mainstream of graph theory. For example, in a proper edge colouring, the subgraphs of the partition are matchings. If each subgraph of a partition is required to be planar (respectively, outerplanar, a forest, a star-forest), then the minimum number of subgraphs in a partition of a graph $G$ is the \emph{thickness} (\emph{outerthickness}, \emph{arboricity}, \emph{star-arboricity}) of $G$. Thickness and arboricity are classical graph parameters that have been studied since the early 1960s. 

The first results in this paper concern the relationship between the above parameters and treewidth, which is a more modern graph parameter that is particularly important in structural and algorithmic graph theory; see the surveys \citep{Bodlaender-TCS98, Reed-AlgoTreeWidth03}. 
In particular, we determine the maximum thickness, maximum outerthickness, maximum arboricity, and maximum star-arboricity of a graph with treewidth $k$. These results are presented in \secref{AbstractResults} (following some background graph theory in \secref{Preliminaries}). 

The main results of the paper are about graph partitions with an additional geometric property. Namely, there is a drawing of the graph, and each subgraph in the partition is drawn without crossings. This type of drawing has  applications in graph visualisation (where each noncrossing subgraph is coloured by a distinct colour), and in multilayer VLSI (where each noncrossing subgraph corresponds to a set of wires that can be routed without crossings in a single layer). With no restriction on how the edges are drawn, the minimum number of noncrossing subgraphs, taken over all drawings of $G$, is again the thickness of $G$. By restricting the edges to be drawn straight, we obtain the \emph{geometric thickness} of $G$. By further restricting the vertices to be in convex position, we obtain the \emph{book thickness} of $G$. These geometric parameters are introduced in \secref{GeometricIntroduction}.

Our main results determine the maximum geometric thickness and maximum book thickness of a graph with treewidth $k$. Analogous results are proved for geometric variations of outerthickness, arboricity, and star-arboricity. These geometric results are stated in \secref{Results}. The general approach that is used in the proofs of our geometric upper bounds is described in \secref{Useful}. The proofs of our geometric results are in \manysecref{BookConstructions}{BookThickness}. \secref{Problems} concludes with numerous open problems.

\section{Background Graph Theory}
\seclabel{Preliminaries}

For undefined graph-theoretic terminology, see the monograph by \citet{Diestel00}. We consider graphs $G$ that are simple, finite, and undirected. Let $V(G)$ and $E(G)$ respectively denote the vertex and edge sets of $G$. For $A,B\subseteq V(G)$, let $G[A;B]$ denote the bipartite subgraph of $G$ with vertex set $A\cup B$ and edge set $\{vw\in E(G):v\in A, w\in B\}$. 

A \emph{partition} of a graph $G$ is a proper partition $\{E_1,E_2,\dots,E_t\}$ of $E(G)$; that is, $\bigcup\{E_i:1\leq i\leq t\}=E(G)$ and $E_i\cap E_j=\emptyset$ whenever $i\ne j$. Each part $E_i$ can be thought of as a spanning subgraph $G_i$ of $G$ with $V(G_i):=V(G)$ and $E(G_i):=E_i$. We also consider a partition to be an edge-colouring, where each edge in $E_i$ is coloured $i$. In an edge-coloured graph, a vertex $v$ is \emph{colourful} if all the edges incident to $v$ receive distinct colours. 

A \emph{graph parameter} is a function $f$ such that $f(G)\in\N:=\{0,1,2,\dots\}$ for every graph $G$. For a graph class $\G$, let $f(\G):=\max\{f(G):G\in\G\}$. If $f(\G)$ is unbounded, we write $f(\G):=\infty$.

Our interest is in drawings of graphs in the plane; see \citep{Szekely-DM04, BrassMoserPach, KaufmannWagner01, DETT99, NishRahman}. A \emph{drawing} $\phi$ of graph $G$ is a pair $(\phi_V,\phi_E)$, where: 

\begin{itemize}
\item $\phi_V$ is an injection from the vertex set $V(G)$ into $\R^2$, and 
\item $\phi_E$ is a mapping from the edge set $E(G)$ into the set of simple curves\footnote{A \emph{simple curve} is a homeomorphic image of the closed unit interval; see \citep{Munkres75} for background topology.} in $\R^2$, such that for each edge $vw\in E(G)$, 
\begin{itemize}
\item the endpoints of the curve $\phi_E(vw)$ are $\phi_V(v)$ and $\phi_V(w)$, and
\item $\phi_V(x)\not\in\phi_E(vw)$ for every vertex $x\in V(G)\setminus\{v,w\}$.
\end{itemize}
\end{itemize}

If $H$ is a subgraph of a graph $G$, then every drawing $\phi$ of $G$ induces a \emph{subdrawing} of $H$ obtained by restricting the functions $\phi_V$ and $\phi_E$ to the elements of $H$. Where there is no confusion, we do not distinguish between a graph element and its image in a drawing. 

A set of points $S\subset\R^2$ is in \emph{general position} if no three points in $S$ are collinear. A \emph{general position} drawing is one in which the vertices are in general position.

Two edges in a drawing \emph{cross} if they intersect at some point other than a common endpoint\footnote{In the literature on crossing numbers it is customary to require that intersecting edges cross `properly' and do not `touch'. This distinction will not be important in this paper.}. A \emph{cell} of a drawing $\phi$ of $G$ is a connected component of $\mathbb{R}^2\setminus\{\phi_V(v):v\in V(G)\}\setminus\cup\{\phi_E(vw):vw\in E(G)\}$. Thus each cell of a drawing is an open subset of $\R^2$ bounded by edges, vertices, and crossing points. Observe that a drawing of a (finite) graph has exactly one cell of infinite measure, called the \emph{outer} cell. A graph drawing with no crossings is \emph{noncrossing}. A graph that admits a noncrossing drawing is \emph{planar}.  A drawing in which all the vertices are on the boundary of the outer cell is \emph{outer}. A graph that admits an outer noncrossing drawing is \emph{outerplanar}.


The \emph{thickness} of a graph $G$, denoted by \tt{G}, is the minimum number of planar subgraphs that partition $G$. Thickness was first defined by \citet{Tutte63a}; see the surveys \citep{Hobbs69, MOS98}. The \emph{outerthickness} of a graph $G$, denoted by $\ot{G}$, is the minimum number of outerplanar subgraphs that partition $G$. Outerthickness was first studied by \citet{Guy74}; also see \citep{Kedlaya-JCTB96, PorMak-CMA05, EC88, Goncalves-STOC05, GN90b, GN90a}. The \emph{arboricity} of a graph $G$, denoted by \arb{G}, is the minimum number of forests that partition $G$. \citet{NW-JLMS64} proved that
\begin{equation}
\eqnlabel{Arboricity}
\arb{G}=\max_{H\subseteq G}\CEIL{\frac{|E(H)|}{|V(H)|-1}}\enspace.
\end{equation}

A \emph{star} is a tree with diameter at most $2$. A \emph{star-forest} is a graph in which each component is a star. The \emph{star-arboricity} of a graph $G$, denoted by \sa{G}, is the minimum number of star-forests that partition $G$. Star arboricity was first studied by \citet{AK85}; also see \cite{Aoki-DM90, AA-DM89, Guiduli-DM97, Huang-PhD, AMR-Comb92, HMS-DM96}.

It is well known that thickness, outerthickness, arboricity, and star-arboricity are within a constant factor of each other. In particular, \citet{Goncalves-STOC05} recently proved a longstanding conjecture that every planar graph $G$ has outerthickness $\ot{G}\leq2$. Thus $\ot{G}\leq2\cdot\tt{G}$ for every graph $G$. Every planar graph $G$ satisfies $|E(G)|<3(|V(G)|-1)$. Thus $\arb{G}\leq3\cdot\tt{G}$ for every graph $G$ by \Eqnref{Arboricity}. Similarly, every outerplanar graph $G$ satisfies $|E(G)|<2(|V(G)|-1)$. Thus $\arb{G}\leq2\cdot\ot{G}$ for every graph $G$ by \Eqnref{Arboricity}. \citet{HMS-DM96} proved that every outerplanar graph $G$ has star-arboricity $\sa{G}\leq3$, and that every planar graph $G$ has star-arboricity $\sa{G}\leq5$. (\citet{AA-DM89} constructed planar graphs $G$ for which $\sa{G}=5$.)\ Thus $\sa{G}\leq3\cdot\ot{G}$  and $\sa{G}\leq5\cdot\tt{G}$ for every graph $G$. It is easily seen that every tree $G$ has star-arboricity $\sa{G}\leq2$. Thus $\sa{G}\leq2\cdot\arb{G}$ for every graph $G$. 
Summarising, we have the following set of inequalities.
\begin{equation}
\eqnlabel{AbstractRelationship}
\begin{minipage}{\textwidth-10mm}
\hspace*{-5mm}
\includegraphics{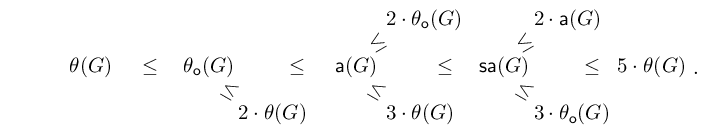}
\end{minipage}
\end{equation}

Let $K_n$ be the complete graph on $n$ vertices. A set of $k$ pairwise adjacent vertices in a graph $G$ is a \emph{$k$-clique}. For a vertex $v$ of $G$, let $N_G(v):=\{w\in V(G):vw\in E(G)\}$ and $N_G[v]:=N_G(v)\cup\{v\}$. We say $v$ is \emph{$k$-simplicial} if $N_G(v)$ is a $k$-clique (and hence $N_G[v]$ is a ($k+1$)-clique). 

For each integer $k\geq1$, a \emph{$k$-tree} is a graph $G$ such that either:\\
\POINT $G\simeq K_{k+1}$, or\\
\POINT $G$ has a $k$-simplicial vertex $v$ and $G\setminus v$ is a $k$-tree.

Suppose that $C$ is a clique in a graph $G$, and $S$ is a nonempty set with $S\cap V(G)=\emptyset$. Let $G'$ be the graph with vertex set $V(G'):=V(G)\cup S$, and edge set $E(G'):=E(G)\cup\{vx:v\in S,x\in C\}$. We say that $G'$ is obtained from $G$ by \emph{adding $S$ onto $C$}. If $S=\{v\}$ then $G'$ is obtained from $G$ by \emph{adding $v$ onto $C$}. Observe that if $|C|=k$, and $G$ is a $k$-tree or $G\simeq K_k$, then $G'$ is a $k$-tree.

The \emph{treewidth} of a graph $G$ is the minimum $k\in\N$ such that $G$ is a spanning subgraph of a $k$-tree. Let $\T_k$ be the class of graphs with treewidth at most $k$. Many families of graphs have bounded treewidth; see \citep{Bodlaender-TCS98}. $\T_1$ is the class of forests. Graphs in $\T_2$ are obviously planar---a $2$-simplicial vertex can always be drawn  near the edge connecting its two neighbours without introducing a crossing. Graphs in $\T_2$ are characterised as those with no $K_4$-minor, and are sometimes called \emph{series-parallel}.

\section{Abstract Parameters and Treewidth}
\seclabel{AbstractResults}




In this section we determine the maximum value of each of thickness, outerthickness, arboricity, and star-arboricity for graphs of treewidth $k$. Since every graph with treewidth $k$ is a subgraph of a $k$-tree, to prove the upper bounds we need only consider $k$-trees. The proofs of the lower bounds employ the \emph{complete split graph} $K^\star_{k,s}$ (for $k,s\geq1$), which is the $k$-tree obtained by adding a set $S$ of $s$ vertices onto an initial $k$-clique $K$; see \figref{Gks}.

\Figure{Gks}{\includegraphics{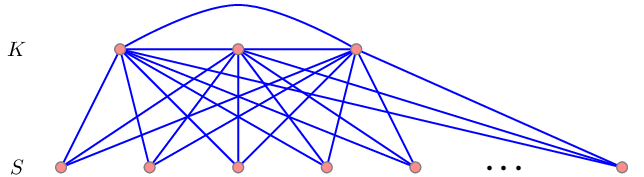}}{The complete split graph $K^\star_{3,s}$.}

Suppose that the edges of $K^\star_{k,s}$ are coloured $1,2,\dots,\ell$. Let $c(e)$ be the colour assigned to each edge $e$ of $K^\star_{k,s}$. The \emph{colour vector} of each vertex $v\in S$ is the set $\{(c(uv),u):u\in K\}$. Note that there are $\ell^k$ possible colour vectors. 

\begin{proposition}
\proplabel{TreewidthThickness}
The maximum thickness of a graph in $\T_k$ is $\ceil{k/2}$; that is, \begin{equation*}\tt{\T_k}=\ceil{k/2}\enspace.
\end{equation*}
\end{proposition}

\begin{proof}
First we prove the upper bound. \citet{DOSV-Comb98} proved that for all $k_1,k_2,\dots,k_t\in\mathbb{N}$ with $k_1+k_2+\dots+k_t=k$, every $k$-tree $G$ can be partitioned into $t$ subgraphs $G_1,G_2,\dots,G_t$, such that each $G_i$ is a $k_i$-tree. Note that the $t=2$ case, which implies the general result, was independently proved by \citet{Chhajed-Networks97}. With $k_i=2$, and since $2$-trees are planar, we have $\tt{G}\leq\ceil{k/2}$. (\thmref{GeometricThickness} provides an alternative proof with additional geometric properties.)

Now we prove the lower bound. If $k\leq2$ then $\tt{\T_k}\geq\tt{K_2}=1=\ceil{k/2}$. Now assume that $k\geq3$. Let $\ell:=\ceil{k/2}-1$ and $s:=2\ell^k+1$. Thus $\ell\geq1$. Suppose that $\tt{K^\star_{k,s}}\leq\ell$. In the corresponding edge $\ell$-colouring of $K^\star_{k,s}$, there are $\ell^k$ possible colour vectors. Thus there are at least three vertices $x,y,z\in S$ with the same colour vector. At least $\ceil{k/\ell}\geq3$ of the $k$ edges incident to $x$ are assigned the same colour. Say these edges are $xa,xb,xc$. Since $y$ and $z$ have the same colour vector as $x$, the $K_{3,3}$ subgraph induced by $\{xa,xb,xc,ya,yb,yc,za,zb,zc\}$ is monochromatic. This is a contradiction since $K_{3,3}$ is not planar. Thus $\tt{\T_k}\geq\tt{K^\star_{k,s}}\geq\ell+1=\ceil{k/2}$.
\end{proof}


\begin{proposition}
\proplabel{TreewidthArboricity}
The maximum arboricity of a graph in $\T_k$ is $k$; that is, 
\begin{equation*}
\arb{\T_k}=k\enspace.
\end{equation*}
\end{proposition}

\begin{proof}
By construction, $|E(G)|=k|V(G)|-k(k+1)/2$ for every $k$-tree $G$. It follows from \Eqnref{Arboricity} that $\arb{G}\leq k$, and $\arb{G}=k$ if $|V(G)|$ is large enough.
\end{proof}


\begin{proposition}
\proplabel{TreewidthOuterThickness}
The maximum outerthickness of a graph in $\T_k$ is $k$; that is, \begin{equation*}\ot{\T_k}=k\enspace.
\end{equation*}
\end{proposition}

\begin{proof}
Since a forest is outerplanar, $\ot{\T_k}\leq\arb{\T_k}=k$ by \propref{TreewidthArboricity}. Now we prove the lower bound. If $k=1$ then $\ot{\T_k}\geq\ot{K_2}=1$. Now assume that $k\geq2$. Let $\ell:=k-1$ and $s:=2\ell^k+1$. Then $\ell\geq1$. Suppose that $\ot{K^\star_{k,s}}\leq\ell$. In the corresponding edge $\ell$-colouring of $K^\star_{k,s}$, there are $\ell^k$ possible colour vectors. Thus there are at least three vertices $x,y,z\in S$ with the same colour vector. At least $\ceil{k/\ell}=2$ of the $k$ edges incident to $x$ are assigned the same colour. Say these edges are $xa$ and $xb$. Since $y$ and $z$ have the same colour vector as $x$, the $K_{2,3}$ subgraph induced by $\{xa,xb,ya,yb,za,zb\}$ is monochromatic. This is a contradiction since $K_{2,3}$ is not outerplanar. Thus $\ot{\T_k}\geq\ot{K^\star_{k,s}}\geq\ell+1=k$.
\end{proof}


\begin{proposition}
\proplabel{TreewidthStarArboricity}
The maximum star-arboricity of a graph in $\T_k$ is $k+1$; that is, \begin{equation*}
\sa{\T_k}=k+1\enspace.
\end{equation*}
\end{proposition}


\begin{proof}
The upper bound $\sa{\T_k}\leq k+1$ was proved by \citet{DOSV-Comb98}\footnote{\lemref{BookStarArboricity} provides an alternative proof that $\sa{\T_k}\leq k+1$. The same result can be concluded from a result by \citet{HMS-DM96}. A vertex colouring with no bichromatic edge and no bichromatic cycle is \emph{acyclic}. It is folklore that every $k$-tree $G$ has an acyclic $(k+1)$-colouring \citep{FRR-JGT04}. (\emph{Proof}. If $G\simeq K_{k+1}$ then the result is trivial. Otherwise, let $v$ be a $k$-simplicial vertex. By induction, $G\setminus v$ has an acyclic $(k+1)$-colouring. One colour is not present on the $k$ neighbours of $v$. Give this colour to $v$. Thus there is no bichromatic edge. The neighbours of $v$ have distinct colours since they form a clique. Thus there is no bichromatic cycle.) \citet{HMS-DM96} proved that a graph with an acyclic $c$-colouring has star arboricity at most $c$. Thus $\sa{\T_k}\leq k+1$.}. For the lower bound, let $s:=k^k+1$. Let $G$ the graph obtained from the $k$-tree $K^\star_{k,s}$ by adding, for each vertex $v\in S$, one new vertex $v'$ onto $\{v\}$. Clearly $G$ has treewidth $k$. Suppose that $\sa{G}\leq k$. In the corresponding edge $k$-colouring of $K^\star_{k,s}$ there are $k^k$ possible colour vectors. Since $|S|>k^k$, there are two vertices $x,y\in S$ with the same colour vector. No two edges in $G[\{x\};K]$ receive the same colour, as otherwise, along with $y$, we would have a monochromatic $4$-cycle. Thus all $k$ colours are present on the edges of $G[\{x\};K]$ and $G[\{y\};K]$. Let $p$ be the vertex in $K$ such that $xp$ and $yp$ receive the same colour as $xx'$. Thus $(x',x,p,y)$ is a monochromatic $4$-vertex path, which is not a star. This contradiction proves that $\sa{\T_k}\geq \sa{G}\geq k+1$. 
\end{proof}


\section{Geometric Parameters}
\seclabel{GeometricIntroduction}

The \emph{thickness} of a graph drawing is the minimum $k\in\N$ such that the edges of the drawing can be partitioned into $k$ noncrossing subdrawings; that is, each edge is assigned one of $k$ colours such that edges with same colour do not cross. Every planar graph can be drawn with its vertices at prespecified locations \citep{Halton91, PW-GC01}. Thus a graph with thickness $k$ has a drawing with thickness $k$ \citep{Halton91}. However, in such a drawing the edges might be highly curved. This motivates the notion of geometric thickness. 

A drawing $(\phi_V,\phi_E)$ of a graph $G$ is \emph{geometric} if the image of each edge $\phi_E(vw)$ is a straight line-segment (by definition, with endpoints $\phi_V(v)$ and $\phi_V(w)$). Thus a geometric drawing of a graph is determined by the positions of its vertices. We thus refer to $\phi_V$ as a geometric drawing.

The \emph{geometric thickness} of a graph $G$, denoted by \gt{G}, is the minimum $k\in\N$ such that there is a geometric drawing of $G$ with thickness $k$. \citet{Kainen73} first defined geometric thickness under the name of \emph{real linear thickness}, and it has also been called \emph{rectilinear thickness}. By the F{\'a}ry-Wagner theorem \citep{Fary48, Wagner36}, a graph has geometric thickness $1$ if and only if it is planar. Graphs of geometric thickness $2$, the so-called \emph{doubly linear} graphs, were studied by \citet{HSV-CGTA99}.

The \emph{outerthickness} (respectively, \emph{arboricity}, \emph{star-arboricity}) of a graph drawing is the minimum $k\in\N$ such that the edges of the drawing can be partitioned into $k$ outer noncrossing subdrawings (noncrossing forests, noncrossing star-forests). Again a graph with outerthickness (arboricity, star-arboricity) $k$ has a drawing with outerthickness (arboricity, star-arboricity) $k$ \citep{Halton91, PW-GC01}. We generalise the notion of geometric thickness as follows. The \emph{geometric outerthickness} (\emph{geometric arboricity}, \emph{geometric star-arboricity}) of a graph $G$, denoted by \got{G} (\ga{G}, \gsa{G}), is the minimum $k\in\N$ such that there is a geometric drawing of $G$ with outerthickness (arboricity, star-arboricity) $k$. 

A geometric drawing in which the vertices are in convex position is called a \emph{book embedding}. The \emph{book thickness} of a graph $G$, denoted by \bt{G}, is the minimum $k\in\N$ such that there is book embedding of $G$ with thickness $k$. Note that whether two edges cross in a book embedding is simply determined by the relative positions of their endpoints in the cyclic order of the vertices around the convex hull. A book embedding with thickness $k$ is commonly called a \emph{$k$-page book embedding}: one can think of the vertices as being ordered on the spine of a book and each noncrossing subgraph being drawn without crossings on a single page. Book embeddings, first defined by \citet{Ollmann73}, are ubiquitous structures with a variety of applications; see \citep{DujWoo-DMTCS04} for a survey with over 50 references. A book embedding is also called a \emph{stack layout}, and book thickness is also called \emph{stacknumber}, \emph{pagenumber} and \emph{fixed outerthickness}.

A graph has book thickness $1$ if and only if it is outerplanar \citep{BK79}. \citet{BK79} proved that a graph has a book thickness at most $2$ if and only if it is a subgraph of a Hamiltonian planar graph. \citet{Yannakakis89} proved that every planar graph has book thickness at most $4$.

The \emph{book arboricity} (respectively, \emph{book star-arboricity}) of a graph $G$, denoted by \ba{G} (\bsa{G}), is the minimum $k\in\N$ such that there is a book embedding of $G$ with arboricity (star-arboricity) $k$. There is no point in defining `book outerthickness' since it would always equal book thickness. By definition,\\[1ex]
\begin{equation*}
\begin{tabular}{ccccc}
\tt{G}	& $\leq$	& \gt{G}	& $\leq$	& \bt{G}\\
$\leq$	&			& $\leq$	& 			& $=$\\
\ot{G}	& $\leq$	& \got{G}	& $\leq$	& \bt{G}\\
$\leq$	&			& $\leq$	& 			& $\leq$\\
\arb{G}	& $\leq$	& \ga{G}	& $\leq$	& \ba{G}\\
$\leq$	&			& $\leq$	& 			& $\leq$\\
\sa{G}	& $\leq$	& \gsa{G}	& $\leq$	& \bsa{G}\enspace.\\
\end{tabular}
\end{equation*}

\section{Main Results}
\seclabel{Results}

As summarised in \tabref{Results}, we determine the value of each geometric graph parameter defined in \secref{GeometricIntroduction} for $\T_k$. 

\begin{table}[h]
\caption{Maximum parameter values for graphs in $\T_k$.}
\tablabel{Results}
\begin{center}
\begin{tabular}{r|cccc}\hline
type of drawing	& thickness		& outerthickness	& arboricity	& star-arboricity	\\\hline
topological		& $\ceil{k/2}$	& $k$				& $k$			& $k+1$				\\			
geometric		& $\ceil{k/2}$	& $k$				& $k$			& $k+1$				\\			
book ($k\leq2$)	& $k$			& -					& $k$			& $k+1$				\\			
book ($k\geq3$)	& $k+1$			& -					& $k+1$			& $k+1$				\\\hline
\end{tabular}
\end{center}
\end{table}

The following theorem is the most significant result in the paper.

\begin{theorem}
\thmlabel{GeometricThickness}
The maximum thickness and maximum geometric thickness of a graph in $\T_k$ satisfy 
\begin{equation*}
\tt{\T_k}=\gt{\T_k}=\ceil{k/2}\enspace.
\end{equation*}
\end{theorem}

\thmref{GeometricThickness} says that the lower bound for the thickness of $\T_k$ (\propref{TreewidthThickness}) can be matched by an upper bound, even in the more restrictive setting of geometric thickness. \thmref{GeometricThickness} is proved in \secref{GeometricDrawings}. 

\begin{theorem}
\thmlabel{GeometricArboricity}
The maximum arboricity, maximum outerthickness, maximum geometric arboricity, and maximum geometric outerthickness of a graph in $\T_k$ satisfy 
\begin{equation*}
\arb{\T_k}=\ot{\T_k}=\got{\T_k}=\ga{\T_k}=k\enspace.
\end{equation*}
\end{theorem}

\thmref{GeometricArboricity} says that our lower bounds for the arboricity and outerthickness of $\T_k$ (\twopropref{TreewidthArboricity}{TreewidthOuterThickness}) can be matched by upper bounds on the corresponding geometric parameter. By the lower bound in \propref{TreewidthOuterThickness}, to prove \thmref{GeometricArboricity}, it suffices to show that $\ga{\T_k}\leq k$; we do so in \secref{GeometricDrawings}.

Now we describe our results for book embeddings. 

\begin{theorem}
\thmlabel{BookThickness}
\thmlabel{BookArboricity}
The maximum book thickness and maximum book arboricity of a graph in $\T_k$ satisfy
\begin{equation*}
\bt{\T_k} = \ba{\T_k} =
\begin{cases}
	k	& \textup{for }k \leq 2\enspace,\\
	k+1 & \textup{for }k \geq 3\enspace.
\end{cases}
\end{equation*}
\end{theorem}

\thmref{BookThickness} with $k=1$ states that every tree has a $1$-page book embedding, as proved by \citet{BK79}. \citet{RM-COCOON95} proved that every series-parallel graph has a $2$-page book embedding (also see \citep{GDLW-Algo06}); that is, $\bt{\T_2}\leq2$. Note that $\bt{\T_2}=2$ since there are series-parallel graphs that are not outerplanar, $K_{2,3}$ being the primary example. We prove the stronger result that $\ba{\T_2}=2$ in \secref{BookConstructions}. 

\citet{GH-DAM01} proved that every $k$-tree has a book embedding with thickness at most $k+1$. In their proof, each noncrossing subgraph is in fact a star-forest. Thus 
\begin{equation}
\eqnlabel{GanleyHeath}
\bt{\T_k}\leq\ba{\T_k}\leq\bsa{\T_k}\leq k+1\enspace.
\end{equation}
We give an alternative proof of this result in \secref{BookConstructions}. \citet{GH-DAM01} proved a lower bound of $\bt{\T_k}\geq k$, and conjectured that $\bt{\T_k}=k$. Thus \thmref{BookThickness} refutes this conjecture. The proof is given in \secref{BookThickness}, where we construct a $k$-tree $Q_k$ with $\bt{Q_k}\geq k+1$. Thus \thmref{BookArboricity} gives an example of an abstract parameter that is not matched by its geometric counterpart. In particular,  $\bt{\T_k}>\ot{\T_k}=k$ for $k\geq3$. 

Note that \citet{TY-DM02} proved that $\bt{G}\leq k$ under the stronger assumption that $G$ has \emph{path}width $k$. Finally observe that the lower bound in \propref{TreewidthStarArboricity} and \Eqnref{GanleyHeath} imply the following result.

\begin{corollary}
\corlabel{StarArboricity}
The maximum star-arboricity, maximum geometric star-arboricity, and maximum book star-arboricity of a graph in $\T_k$ satisfy
\begin{equation*}\sa{\T_k}=\gsa{\T_k}=\bsa{\T_k}=k+1\enspace.
\end{equation*}
\end{corollary}

\section{General Approach}
\seclabel{Useful}

When proving upper bounds, we need only consider $k$-trees, since edges can be added to a graph with treewidth $k$ to obtain a $k$-tree, without decreasing the relevant thickness or arboricity parameter. The definition of a $k$-tree $G$ suggests a natural approach to drawing $G$: choose a simplicial vertex $w$, recursively draw $G\setminus w$, and then add $w$ to the drawing. For the problems under consideration this approach fails because the neighbours of $w$ may have high degree. The following lemma solves this impasse.

\begin{lemma}
\lemlabel{Partition}
Every $k$-tree $G$ has a nonempty independent set $S$ of $k$-simplicial vertices such that either:
\begin{enumerate}
\item[\textup{(a)}] $G\setminus S\simeq K_k$ \paran{that is, $G\simeq K^\star_{k,|S|}$}, or
\item[\textup{(b)}] $G\setminus S$ is a $k$-tree containing a $k$-simplicial vertex $v$ such that:
\begin{itemize}
\item for each vertex $w\in S$, there is exactly one vertex $u\in N_{G\setminus S}(v)$ such that $N_G(w)=N_{G\setminus S}[v]\setminus\{u\}$, and
\item each $k$-simplicial vertex of $G$ that is not in $S$ is not adjacent to $v$.
\end{itemize}
\end{enumerate}
\end{lemma}

\begin{proof}
Every $k$-tree has at least $k+1$ vertices. If $|V(G)|=k+1$ then $G\simeq K_{k+1}$ and property (a) is satisfied with $S=\{v\}$ for each vertex $v$. Now assume that $|V(G)|\geq k+2$. Let $L$ be the set of $k$-simplicial vertices of $G$. Then $L$ is a nonempty independent set, and $G\setminus L$ is a $k$-tree or $G\setminus L\simeq K_k$. If $G\setminus L\simeq K_k$, then property (a) is satisfied with $S=L$. Otherwise, $G\setminus L$ has a $k$-simplicial vertex $v$. Let $S$ be the set of neighbours of $v$ in $L$. We claim that property (b) is satisfied. Now $S\ne\emptyset$, as otherwise $v\in L$. Since $G$ is not a clique and each vertex in $S$ is simplicial, $G\setminus S$ is a $k$-tree. Consider a vertex $w\in S$. Now $N_G(w)$ is a $k$-clique and $v\in N_G(w)$. Thus $N_G(w)\subseteq N_{G\setminus S}[v]$. Since $|N_G(w)|=k$ and $|N_{G\setminus S}[v]|=k+1$, there is exactly one vertex $u\in N_{G\setminus S}(v)$ for which $N_G(w)=N_{G\setminus S}[v]\setminus\{u\}$. The final claim is immediate.
\end{proof}

\lemref{Partition} is used to prove all of the upper bounds that follow. Our general approach is: 
\begin{itemize}
\item in a recursively computed drawing of $G\setminus S$, draw the vertices in $S$ close to $v$,
\item for each vertex $w\in S$, colour the edge $wx$ ($x\ne v$) by the colour assigned to $vx$, and colour the edge $wv$ by the colour assigned to the edge $vu$, where $u$ is the neighbour of $v$ that is not adjacent to $w$.
\end{itemize}

\section{Constructions of Book Embeddings}
\seclabel{BookConstructions}

First we prove that $\bsa{\T_k}=k+1$. The lower bound follows from the stronger lower bound $\sa{\T_k}\geq k+1$ in \propref{TreewidthStarArboricity}. The upper bound is proved by induction on $|V(G)|$ with the following hypothesis. Recall that in an edge-coloured graph, a vertex $v$ is \emph{colourful} if all the edges incident to $v$ receive distinct colours. 

\begin{lemma}
\lemlabel{BookStarArboricity}
Every $k$-tree $G$ has a book embedding with star-arboricity $k+1$ such that:\\
\POINT if $G\simeq K_{k+1}$ then at least one vertex is colourful, and \\
\POINT if $G\not\simeq K_{k+1}$ then every $k$-simplicial vertex is colourful.
\end{lemma}

\begin{proof}
Apply \lemref{Partition} to $G$. We obtain a nonempty independent set $S$ of $k$-simplicial vertices of $G$. 

First suppose that $G\setminus S\simeq K_k$ with $V(G\setminus S)=\{u_1,u_2,\dots,u_k\}$. Position $V(G)$ arbitrarily on a circle, and draw the edges straight. Every edge of $G$ is incident to some $u_i$. Colour the edges $1,2,\dots,k$ so that every edge coloured $i$ is incident to $u_i$. Thus each colour class is a noncrossing star, and every vertex in $S$ is colourful. If $G\simeq K_{k+1}$ then $|S|=1$ and at least one vertex is colourful. If $G\not\simeq K_{k+1}$ then no vertex $u_i$ is $k$-simplicial; thus every $k$-simplicial vertex is in $S$ and is colourful.

Otherwise, by \lemref{Partition}(b), $G\setminus S$ is a $k$-tree containing a $k$-simplicial vertex $v$, such that $N_G(w)\subset N_{G\setminus S}[v]$ for each vertex $w\in S$. Say $N_{G\setminus S}(v)=\{u_1,u_2,\dots,u_k\}$. 

Apply the induction hypothesis to $G\setminus S$. If $G\setminus S\simeq K_{k+1}$ then we can nominate $v$ to be a vertex of $G\setminus S$ that becomes colourful. By induction, we obtain a book embedding of $G\setminus S$ with star-arboricity $k+1$, in which $v$ is colourful. Without loss of generality, each edge $vu_i$ is coloured $i$.  Let $x$ be a vertex next to $v$ on the convex hull. Position the vertices in $S$ arbitrarily between $v$ and $x$. For each $w\in S$, colour each edge $wu_i$ by $i$, and colour $wv$ by $k+1$, as illustrated in \figref{BookConstructions}(a).

\Figure{BookConstructions}{\includegraphics{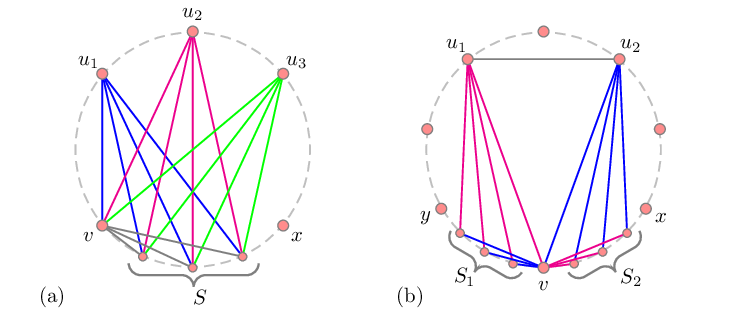}}{Book embedding constructions: (a) in \lemref{BookStarArboricity} with $k=3$, and (b) in \lemref{BookArboricity}.}

By construction, each vertex in $S$ is colourful. The edges $\{vw:w\in S\}$ form a new star component of the star-forest coloured $k+1$. For each colour $i\in\{1,2,\dots,k\}$, the component of the subgraph of $G\setminus S$ that is coloured $i$ and contains $v$ is a star rooted at $u_i$ with $v$ a leaf. Thus it remains a star by adding the edge $wu_i$ for all $w\in S$. 

Suppose that two edges $e$ and $f$ of $G$ cross and are both coloured $i$ ($\in\{1,2,\dots,k\}$). Then $e$ and $f$ are not both in $G\setminus S$. Without loss of generality, $e$ is incident to a vertex $w\in S$. The edges of $G$ that are coloured $i$ and have at least one endpoint in $S\cup\{v\}$ form a noncrossing star (rooted at $u_i$ if $1\leq i\leq k$, and rooted at $v$ if $i=k+1$). Thus $f$ has no endpoint in $S\cup\{v\}$. Observe that $vw$ crosses no edge in $G\setminus S$. Thus $e=wu_i$. Since $S\cup\{v\}$ is consecutive on the circle and $f$ has no endpoint in $S\cup\{v\}$, $f$ also crosses $vu_i$. Hence $f$ and $vu_i$ are two edges of $G\setminus S$ that cross and are both coloured $i$. This contradiction proves that no two edges of $G$ cross and receive the same colour.

It remains to prove that every $k$-simplicial vertex in $G$ is colourful.
Each vertex in $S$ is colourful. Consider a $k$-simplicial vertex $x$ of $G$ that is not in $S$. By \lemref{Partition}(b), $x$ is not adjacent to $v$. Thus $x$ is adjacent to no vertex in $S$, and $x$ is $k$-simplicial in $G\setminus S$. Moreover, $G\setminus S$ is not complete. By induction, $x$ is colourful in $G\setminus S$ and in $G$.
\end{proof}


Now we prove \thmref{BookArboricity} with $k=2$, which states that $\bt{\T_2}=\ba{\T_2}=2$. The lower bound $\bt{\T_2}\geq2$ holds since $K_{2,3}$ is series-parallel but is not outerplanar. We prove the upper bound $\ba{\T_2}\leq2$ by induction on $|V(G)|$ with the following hypothesis.

\begin{lemma}
\lemlabel{BookArboricity} 
Every $2$-tree $G$ has a book embedding with arboricity $2$ such that:\\
\POINT if $G\simeq K_3$ then two vertices are colourful, and\\
\POINT if $G\not\simeq K_3$ then every $2$-simplicial vertex is colourful.
\end{lemma}

\begin{proof}
Apply \lemref{Partition} to $G$. We obtain a nonempty independent set $S$ of $2$-simplicial vertices of $G$. 

First suppose that $G\setminus S\simeq K_2$ with $V(G\setminus S)=\{u_1,u_2\}$. 
Position $V(G)$ at distinct points on a circle in the plane, and draw the edges straight. Every edge is incident to $u_1$ or $u_2$. Colour every edge incident to $u_1$ by $1$. Colour every edge incident to $u_2$ (except $u_1u_2$) by $2$. Thus each colour class is a noncrossing star, and each vertex in $S$ is colourful. If $G\simeq K_3$ then $|S|=1$ and $u_2$ is also colourful. If $G\not\simeq K_3$ then neither $u_1$ nor $u_2$ are $2$-simplicial; thus each $2$-simplicial vertex is colourful.

Otherwise, by \lemref{Partition}(b), $G\setminus S$ is a $2$-tree containing a $2$-simplicial vertex $v$. Say $N_{G\setminus S}(v)=\{u_1,u_2\}$. For every vertex $w\in S$, $N_G(w)=\{v,u_1\}$ or $N_G(w)=\{v,u_2\}$. Let $S_1=\{w\in S:N_G(w)=\{v,u_1\}\}$ and  $S_2=\{w\in S:N_G(w)=\{v,u_2\}\}$. 

Apply the induction hypothesis to $G\setminus S$. If $G\setminus S\simeq K_3$ we can nominate $v$ to be a vertex of $G\setminus S$ that becomes colourful. By induction, we obtain a book embedding of $G\setminus S$ with arboricity $2$, in which $v$ is colourful. Without loss of generality, each edge $vu_i$ is coloured $i$. Say $u_1$ appears before $u_2$ in clockwise order from $v$. Say $(x,v,y)$ are consecutive in clockwise order, as illustrated in \figref{BookConstructions}(b). Position the vertices in $S_1$ between $v$ and $y$, and position the vertices in $S_2$ between $x$ and $v$. For all $w\in S_i$, colour each edge $wu_i$ by $i$, and colour $wv$ by $3-i$. 

The only edge that can cross an edge $wv$ ($w\in S_i$) is some $pu_i$ where $p\in S_i$. These edges receive distinct colours. If an edge $e$ of $G\setminus S$ crosses some edge $wu_i$, then $e$ also crosses $vu_i$ (since $\deg_{G\setminus S}(v)=2$). Since $wu_i$ receives the same colour as $vu_i$, $e$ must be coloured differently from $wu_i$. Hence edges assigned the same colour do not cross. 

By construction, each vertex $w\in S$ is colourful; $w$ becomes a leaf in both forests of the partition. It remains to prove that every $2$-simplicial vertex in $G$ is colourful. Each vertex in $S$ is colourful. Consider a $k$-simplicial vertex $x$ of $G$ that is not in $S$. By \lemref{Partition}(b), $x$ is not adjacent to $v$. Thus $x$ is adjacent to no vertex in $S$, and $x$ is $2$-simplicial in $G\setminus S$. Moreover, $G\setminus S$ is not complete. By induction, $x$ is colourful in $G\setminus S$ and in $G$.
\end{proof}

\section{Constructions of Geometric Drawings}
\seclabel{GeometricDrawings}

In this section we prove \twothmref{GeometricThickness}{GeometricArboricity}. First we introduce some geometric notation. Let $v$ and $w$ be distinct points in the plane; see \figref{Lines}. 
Let \lin{v}{w} be the line through $v$ and $w$. 
Let \seg{v}{w} be the open line-segment with endpoints $v$ and $w$. 
Let \segc{v}{w} be the closed line-segment with endpoints $v$ and $w$.
Let \ray{v}{w} be the open ray from $v$ through $w$. 
Let \oray{v}{w} be the open ray opposite to \ray{v}{w}; that is,  $\oray{v}{w}:=(\lin{v}{w}\setminus\ray{v}{w})\setminus\{v\}$.
 
\Figure{Lines}{\includegraphics{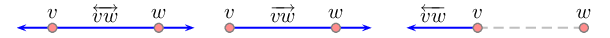}}{Notation for lines and rays.}

For every point $p\in\R^2$ and set of points $Q\subset\R^2\setminus\{p\}$, such that $Q\cup\{p\}$ is in general position, let 
\begin{equation*}
R(p,Q):=\{\ray{p}{q},\oray{p}{q}:q\in Q\}
\end{equation*} 
be the set of rays from $p$ to the points in $Q$ together with their opposite rays, in clockwise order around $p$. (Since $Q\cup\{p\}$ is in general position, the rays in $R(p,Q)$ are pairwise disjoint, and their clockwise order is unique.)\ 

Let $r$ and $r'$ be non-collinear rays from a single point $v$. The \emph{wedge} \WEDGE{r}{r'} \emph{centred} at $v$ is the unbounded region of the plane obtained by sweeping a ray from $r$ to $r'$ through the lesser of the two angles formed by $r$ and $r'$ at $v$. We consider \WEDGE{r}{r'} to be open in the sense that $r\cup r'\cup\{v\}$ does not intersect \WEDGE{r}{r'}.

The proofs of \twothmref{GeometricThickness}{GeometricArboricity} are incremental constructions of geometric drawings. 
The insertion of new vertices is based on the following definitions. 

Consider a geometric drawing of a graph $G$. 
Let $v$ be a vertex of $G$. 
For $\varepsilon>0$, let $D_\varepsilon(v)$ be the open disc of radius $\varepsilon$ centred at $v$. For a point $u$, let 
\begin{equation*}
C_\varepsilon(v,u):=\bigcup\{\seg{u}{x}:x\in D_\varepsilon(v)\}
\end{equation*} 
be the region in the plane consisting of the union of all open line-segments from $u$ to the points in $D_\varepsilon(v)$. Let 
\begin{equation*}
T_\varepsilon(v):=\bigcup\{C_\varepsilon(v,u):u\in N_G(v)\}
\end{equation*} 
be the region in the plane consisting of the union of all open line-segments from each neighbour of $v$ to the points in $D_\varepsilon(v)$. 

As illustrated in \figref{Construction}(a), a vertex $v$ in a general position geometric drawing of a graph $G$ is \emph{$\varepsilon$-empty} if:\\
\aaa\ the only vertex of $G$ in $T_\varepsilon(v)$ is $v$,\\
\bbb\ every edge of $G$ that intersects $D_\varepsilon(v)$ is incident to $v$,\\
\ccc\ $(V(G)\setminus\{v\})\cup\{p\}$ is in general position for each point $p\in D_\varepsilon(v)$, and \\
\ddd\ the clockwise orders of $R(v,N_G(v))$ and $R(p,N_G(v))$ are the same\\ 
\hspace*{1.8em} for each point $p\in D_\varepsilon(v)$.

\Figure{Construction}{\includegraphics{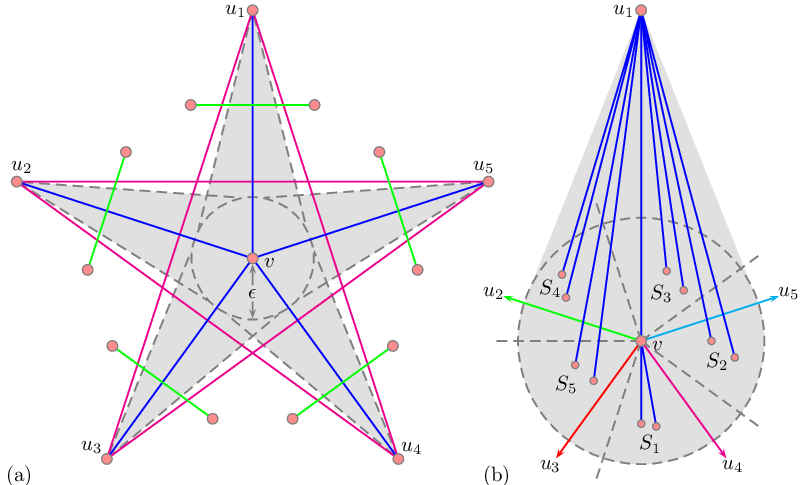}}{(a) $N_G(v)=\{u_1,u_2,\dots,u_5\}$, $T_\varepsilon(v)$ is shaded, and $v$ is $\varepsilon$-empty. (b) The new edges coloured $1$ in \propref{GeometricArboricity}.}


We have the following observations.

\begin{observation}
\obslabel{EmptyDisc}
Every vertex $v$ in a general position geometric drawing of a graph $G$ is $\varepsilon$-empty for some $\varepsilon>0$.
\end{observation}

\begin{proof}
Consider the arrangement $A$ consisting of the lines through every pair of vertices in $G\setminus v$; see \citep{Mat02} for background on line arrangements. Since $V(G)$ is in general position, $v$ is in some cell $C$ of $A$. Since $C$ is an open set, there exists $\varepsilon>0$ such that $D_\varepsilon(v)\subset C$. For every neighbour $u\in N_G(v)$, no vertex $x$ of $G\setminus v$ is in $C_\varepsilon(v,u)$, as otherwise \lin{x}{u} would intersect $D_\varepsilon(v)$. Thus property (a) holds. No line of $A$ intersects $C$. In particular, no edge of $G\setminus v$ intersects $C$, and property (b) holds. No point $p\in D_\varepsilon(v)$ is collinear with two vertices of $G\setminus v$, as otherwise $D_\varepsilon(v)$ would intersect a line in $A$. Thus property (c) holds. The radial order of $V(G)\setminus v$ is the same from each point in $C$. In particular, property (d) holds. Therefore $v$ is $\varepsilon$-empty.
\end{proof}


\begin{observation}
\obslabel{MustCross}
Let $v$ be an $\varepsilon$-empty vertex in a general position geometric drawing of a graph $G$. Let $u\in N_G(v)$. Suppose that some edge $e\in E(G)$ crosses \seg{p}{u} for some point $p\in D_\varepsilon(v)$. Then either $e$ is incident to $v$, or $e$ also crosses the edge $vu$.
\end{observation}

\begin{proof}
If $e$ is incident to $v$, then we are done. Now assume that $e$ is not incident to $v$. Thus $e$ does not intersect $D_\varepsilon(v)$ by property (b) of the choice of $\varepsilon$. Since $p\in D_\varepsilon(v)$, we have $\seg{p}{u}\subset C_\varepsilon(v,u)$. Thus the crossing point between $e$ and $\seg{p}{u}$ is in $C_\varepsilon(v,u)\setminus D_\varepsilon(v)$. In particular, $e$ intersects $C_\varepsilon(v,u)$. By property (a) of the choice of $\varepsilon$ and since $e$ is not incident to $v$, no endpoint of $e$ is in $T_\varepsilon(v)$.

We have proved that $e$ does not intersect $D_\varepsilon(v)$, $e$ intersects $C_\varepsilon(v,u)$, and no endpoint of $e$ is in $T_\varepsilon(v)$. Observe that any segment with these three properties must cross $vu$. Thus $e$ crosses $vu$.
\end{proof}

\subsection{Proof of \thmref{GeometricArboricity}}




\thmref{GeometricArboricity} states that $\arb{\T_k}=\ot{\T_k}=\got{\T_k}=\ga{\T_k}=k$. By the discussion in \secref{Results}, it suffices to show that for geometric arboricity,  $\ga{\T_k}\leq k$. We proceed by induction on $|V(G)|$ with the following hypothesis. 

\begin{proposition} 
\proplabel{GeometricArboricity}
Every $k$-tree $G$ has a general position geometric drawing with arboricity $k$ such that:\\
\POINT if $G\simeq K_{k+1}$ then at least one vertex is colourful, and \\
\POINT if $G\not\simeq K_{k+1}$ then every $k$-simplicial vertex is colourful.
\end{proposition}

\begin{proof}
Apply \lemref{Partition} to $G$. We obtain a nonempty independent set $S$ of $k$-simplicial vertices of $G$. 

First suppose that $G\setminus S\simeq K_k$ with $V(G\setminus S)=\{u_1,u_2,\dots,u_k\}$. Fix an arbitrary general position geometric drawing of $G$. Greedily colour the edges of $G$ with colours $1,2,\dots,k$, starting with the edges incident to $u_1$ and ending with the edges incident to $u_k$, so that every edge coloured $i$ is incident to $u_i$. Thus each colour class is a noncrossing star, and every vertex in $S$ is colourful. If $G\simeq K_{k+1}$ then $|S|=1$ and at least one vertex is colourful. If $G\not\simeq K_{k+1}$ then no vertex $u_i$ is $k$-simplicial in $G$; thus each $k$-simplicial vertex is in $S$ and is colourful.

Otherwise, by \lemref{Partition}(b), $G\setminus S$ is a $k$-tree containing a $k$-simplicial vertex $v$. Say $N_{G\setminus S}(v)=\{u_1,u_2,\dots,u_k\}$. Each vertex $w\in S$ has $N_G(w)=N_{G\setminus S}[v]\setminus\{u_i\}$ for exactly one value of $i\in\{1,2,\dots,k\}$. Let $S_i:=\{w\in S:N_G(w)=N_{G\setminus S}[v]\setminus\{u_i\}\}$ for each $i\in\{1,2,\dots,k\}$. Then $\{S_1,S_2,\dots,S_k\}$ is a partition of $S$.

Apply the induction hypothesis to $G\setminus S$. If $G\setminus S\simeq K_{k+1}$ then we can nominate $v$ to be a vertex of $G\setminus S$ that becomes colourful. By induction, we obtain a general position geometric drawing of $G\setminus S$ with arboricity $k$, in which $v$ is colourful. Without loss of generality, each edge $vu_i$ is coloured $i$. 

By \obsref{EmptyDisc}, $v$ is $\varepsilon$-empty in the general position geometric drawing of $G\setminus S$ for some $\varepsilon>0$. Let $X_1,X_2,\dots,X_k$ be pairwise disjoint wedges centred at $v$ such that 
$\oray{v}{u_i}\subset X_i$ for all $i\in\{1,2,\dots,k\}$. Position the vertices of $S_i$ in $X_i\cap D_\varepsilon(v)$ so that $V(G)$ is in general position. This is possible since $X_i\cap D_\varepsilon(v)$ is an open (infinite) region, but there are only finitely many pairs of vertices. Draw each edge straight. For each vertex $w\in S_i$, colour the edge $wv$ by $i$, and colour the edge $wu_j$ ($j\ne i$) by $j$. Thus $w$ is colourful; $w$ becomes a leaf in each of the $k$ forests. This construction is illustrated in \figref{Construction}(b).

To prove that edges assigned the same colour do not cross, consider the set of edges coloured $i$ to be partitioned into three sets: \\
\ONE\ edges in $G\setminus S$ that are coloured $i$,\\
\TWO\ edges $wu_i$ for some $w\in S\setminus S_i$, and \\
\THREE\ edges $vw$ for some $w\in S_i$. 

Type-(1) edges do not cross by induction. Type-(2) edges do not cross since they are all incident to $u_i$. Type-(3) edges do not cross since they are all incident to $v$. 

Suppose that a type-(1) edge $e$ crosses a type-(2) edge $wu_i$ for some $w\in S$. By \obsref{MustCross} with $p=w$ ($\in D_\varepsilon(v)$), either $e$ is incident to $v$, or $e$ also crosses $vu_i$. Since $e$ and $vu_i$ are both coloured $i$, they do not cross in $G$, and we can now assume that $e$ is incident to $v$. Thus $e=vu_i$, which is the only edge in $G\setminus S$ that is incident to $v$ and is coloured $i$. Since $e$ and $wu_i$ have a common endpoint, $e$ and $wu_i$ do not cross, which is a contradiction. Thus type-(1) and type-(2) edges do not cross. 

Now suppose that a type-(1) edge $e$ crosses a type-(3) edge $wv$ for some $w\in S_i$. Then $e\ne vu_i$, since $vu_i$ and $wv$ have a common endpoint. Now, $wv$ is contained in $D_\varepsilon(v)$. Thus $e$ intersects $D_\varepsilon(v)$, which contradicts property (b) of the choice of $\varepsilon$. Thus type-(1) and type-(3) edges do not cross. 

By construction, no type-(2) edge intersects the wedge $X_i$. Since every type-(3) edge is contained in $X_i$, type-(2) and type-(3) edges do not cross. Therefore edges assigned the same colour do not cross. 

It remains to prove that each $k$-simplicial vertex of $G$ is colourful. Each vertex in $S$ is colourful. Consider a $k$-simplicial vertex $x$ that is not in $S$. By \lemref{Partition}(b), $x$ is not adjacent to $v$. Thus $x$ is adjacent to no vertex in $S$, and $x$ is $k$-simplicial in $G\setminus S$. Moreover, $G\setminus S$ is not complete. By induction, $x$ is colourful in $G\setminus S$ and in $G$.
\end{proof}

\subsection{Proof of \thmref{GeometricThickness}}



\thmref{GeometricThickness} states that $\tt{\T_k}=\gt{\T_k}=\ceil{k/2}$. The thickness lower bound, $\tt{\T_k}\geq\ceil{k/2}$, is \propref{TreewidthThickness}. For the the upper bound on the geometric thickness, $\gt{\T_k}\leq\ceil{k/2}$, it suffices to prove that $\gt{\T_{2k}}\leq k$ for all $k\geq2$ (since graphs in $\T_2$ are planar, and thus have geometric thickness $1$). We use the following definitions, for some fixed $k\geq 2$. Let 
\begin{equation*}
I:=\{\ii,-\ii:1\leq \ii\leq k\}\enspace.
\end{equation*} 

Suppose that $\phi$ is a geometric drawing of a graph $G$. (Note that $G$ is not necessarily a $2k$-tree, and $\phi$ is not necessarily in general position.)\ Suppose that $v$ is a vertex of $G$ with degree $2k$, where 
\begin{equation}
\eqnlabel{NeighbourList}
N_G(v)=(u_1,u_2,\dots,u_k,u_{-1},u_{-2},\dots,u_{-k})
\end{equation} 
are the neighbours of $v$ in clockwise order around $v$ in $\phi$. (Since no edge passes through a vertex, this cyclic ordering is well defined.)\ For each $\ii\in I$, define the \emph{$i$-wedge} of $v$ (with respect to the labelling of $N_G(v)$ in \Eqnref{NeighbourList}) to be
\begin{equation*}
F_\ii(v):=\WEDGE{\ray{v}{u_\ii}}{\oray{v}{u_{-\ii}}}\enspace.
\end{equation*} 
If $u_\ii,v,u_\jj$ are collinear, then $\ray{v}{u_\ii}=\oray{v}{u_\jj}$. But if $\phi$ is in general position, then $\ray{v}{u_\ii}\ne\oray{v}{u_\jj}$ for all $\ii,\jj\in I$. Now suppose that, in addition, $\phi$ is in general position. Let
\begin{equation*}
R(v):=R(v,N_G(v))=\{\ray{v}{u_\ii},\oray{v}{u_\ii}:\ii\in I\}
\end{equation*} 
be the set of $2k$ open rays from $v$ through its neighbours together with their $2k$ opposite open rays, in clockwise order around $v$ in $\phi$. We say $v$ is \emph{balanced} in $\phi$ if \ray{v}{u_\ii} and \oray{v}{u_{-\ii}} are consecutive in $R(v)$ for each $\ii\in I$. Note that $v$ is balanced if and only if $F_\ii(v)\cap F_\jj(v)=\emptyset$ for all distinct $\ii,\jj\in I$. Moreover, whether $v$ is balanced does not depend on the choice of labelling in \Eqnref{NeighbourList}.

Now suppose that, in addition, $G$ is a $2k$-tree, and $\phi$ has thickness $k$. Consider the edges of $G$ to be coloured $1,2,\dots,k$, where edges of the same colour do not cross in $\phi$. As illustrated in \figref{ShadowsDuplicate}(a), a $2k$-simplicial vertex $v$ of $G$ is a \emph{fan} in $\phi$ if, for some labelling of $N_G(v)$ as in \Eqnref{NeighbourList}, we have:\\
\POINT $v$ is balanced in $\phi$, and  \\
\POINT the edge $vu_\ii$ is coloured $|\ii|$ for each $\ii \in I$.

Note that for all $Q\subseteq V(G)$ and all $v\in Q$ such that $G[Q]$ is a $2k$-tree and $v$ is $2k$-simplicial\footnote{Since $G[Q]$ is a $2k$-tree, it has minimum degree $2k$. Since $v\in Q$ and $\deg_G(v)=2k$, we have $\deg_{G[Q]}(v)=2k$. Thus every neighbour of $v$ in $G$ is also in $Q$. Thus $v$ is $2k$-simplicial in $G[Q]$.} in $G$, $v$ is a fan in $\phi$ if and only if $v$ is a fan in the drawing of $G[Q]$ induced by $\phi$.

A drawing $\phi$ of a $2k$-tree $G$ is \emph{good} if:\\
\POINT $\phi$ is a general position geometric drawing, \\
\POINT $\phi$ has thickness $k$, \\ 
\POINT if $G\simeq K_{2k+1}$ then at least one vertex of $G$ is a fan in $\phi$, and\\
\POINT if $G\not\simeq K_{2k+1}$ then every $2k$-simplicial vertex of $G$ is a fan in $\phi$.

\Figure{ShadowsDuplicate}{\includegraphics{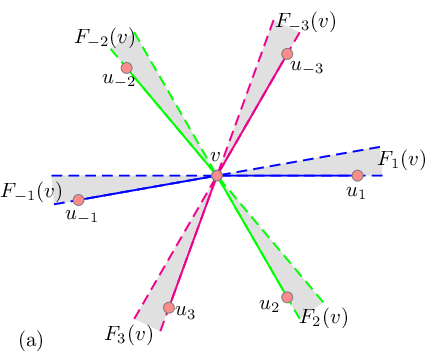}\includegraphics{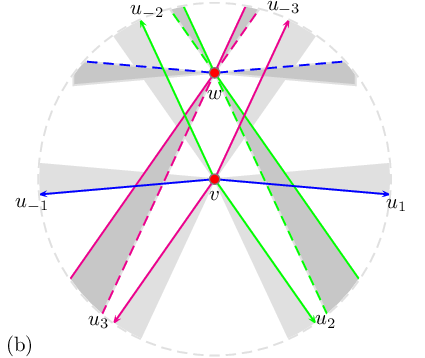}}{(a) A fan vertex $v$ with $k=3$. (b) Inserting the vertex $w$ in \lemref{Insert}.}

The proof of \thmref{GeometricThickness} uses the following two lemmas about constructing good drawings.

\begin{lemma}
\lemlabel{Insert}
Consider a $2k$-tree $G$ for some $k\geq 2$. Suppose that $G$ has a good drawing $\phi$, and $v$ is a fan vertex in $\phi$. Let $G'$ be the $2k$-tree obtained from $G$ by adding a new vertex $w$ onto $N_G(v)$. Then $w$ can be inserted into $\phi$ to obtain a good drawing $\phi'$ of $G'$.
\end{lemma}

\begin{proof}
Say $(u_1,u_2,\dots,u_k,u_{-1},u_{-2},\dots,u_{-k})$ are the neighbours of $v$ in clockwise order around $v$. Since $v$ is a fan in $\phi$, the edge $vu_\ii$ is coloured $|\ii|$ for all $\ii\in I$. By \obsref{EmptyDisc}, $v$ is $\varepsilon$-empty for some $\varepsilon>0$. 
Let 
\begin{equation*}
X:=D_\varepsilon(v)\setminus\{v\}\setminus\cup\{F_\ii(v):\ii\in I\}\enspace.
\end{equation*} 
Thus $X$ consists of $2k$ connected sets having nonempty interior. Hence, there is a nonempty, in fact open, subset of $X$ consisting of points that are not collinear with any two distinct vertices of $G$. Map $w$ to any point in that subset, and draw each edge $wu_\ii$ straight ($\ii\in I$). We obtain a general position geometric drawing $\phi'$ of $G'$.
As illustrated in \figref{ShadowsDuplicate}(b), colour each edge $wu_\ii$ of $G'$ by $|\ii|$, which is the same colour assigned to $vu_\ii$. 

Consider an edge $e$ of $G$ that crosses $wu_\ii$ in $\phi'$ for some $\ii\in I$. By construction, $wu_\ii$ is coloured $|\ii|$. Suppose, for the sake of contradiction, that $e$ is also coloured $|\ii|$. By \obsref{MustCross} with $p=w$ ($\in D_\varepsilon(v)$), either $e$ is incident to $v$, or $e$ also crosses $vu_\ii$. Since $e$ and $vu_\ii$ are both coloured $i$ in $G$, $e$ does not cross $vu_\ii$, and we can now assume that $e$ is incident to $v$. Since $vu_\ii$ and $wu_\ii$ share an endpoint, $e\ne vu_\ii$. Thus $e=vu_{-\ii}$, which is the only other edge incident to $v$ coloured $|\ii|$. Since $wu_\ii$ crosses $vu_{-\ii}$, we have that $w\in F_{-\ii}(v)$, which contradicts the placement of $w$. Thus edges of $G'$ that are assigned the same colour do not cross in $\phi'$. 

Let $x\ne w$ be a $2k$-simplicial vertex in $G'$. Then $x$ is not adjacent to $w$, and $x$ is $2k$-simplicial in $G$. Since $x$ is a fan in $\phi$, it also is a fan in $\phi'$. We now prove that $w$ is a fan in $\phi'$. By property (d) of the choice of $\varepsilon$, and since $w\in D_\varepsilon(v)$, the cyclic orderings of the ray sets $R(v)$ and $R(w)$ are the same. Since $v$ is a fan in $\phi$, and by the colouring of the edges incident to $w$, $w$ is also fan in $\phi'$. 

If $G\simeq K_{2k+1}$, then $v$ and $w$ are the only $2k$-simplicial vertices in $G'$, and thus every $2k$-simplicial  vertex of $G'$ is a fan in $\phi'$. If $G\not\simeq K_{2k+1}$, consider a $2k$-simplicial vertex $y\not=w$ of $G'$. No pair of $2k$ simplicial vertices in $G'$ are adjacent. Thus $y$ is $2k$-simplicial in $G$ and $y$ is a fan in $\phi$ (and $\phi'$). Thus every $2k$-simplicial vertex of $G'$ is a fan in $\phi'$, as required.
\end{proof}

\begin{lemma}
\lemlabel{BaseCase}
For all $k\geq2$, the complete graph $K_{2k+1}$ has a good drawing in which any given vertex $v$ is a fan.
\end{lemma}

\begin{proof}
Say $V(K_{2k+1})=\{v,u_1,u_2,\dots,u_{2k}\}$. As illustrated in \figref{BaseCase}(a), position $u_1,u_2,\dots,u_{2k}$ evenly spaced and in this order on a circle in the plane centred at a point $p$. The edges induced by $\{u_1,u_2,\dots,u_{2k}\}$ can be $k$-coloured using the standard book embedding of $K_{2k}$ with thickness $k$: colour each edge $u_\alpha u_\beta$ by $1+\floor{\half((\alpha+\beta)\bmod{2k})}$. Then the colours are $1,2,\dots,k$, and each colour class forms a noncrossing zig-zag subgraph \citep{BK79,DujWoo-DMTCS04}.

Rename each vertex $u_{k+i}$ by $u_{-\ii}$. As illustrated in \figref{BaseCase}(b), the edges $\{u_\ii u_{-\ii}:1\leq \ii\leq k\}$ pairwise intersect at $p$. Position $v$ strictly inside a cell of the drawing of $K_{2k}$ that borders $p$ (the shaded region in \figref{BaseCase}(a)). Then $V(K_{2k+1})$ is in general position. For all $\ii\in I$, colour $vu_\ii$ by $|\ii|$. Then edges assigned the same colour do not cross. $v$ is a fan since $R(v,\{u_\ii:\ii\in I\})=(\oray{v}{u_{-1}},\ray{v}{u_1},\oray{v}{u_{-2}},\ray{v}{u_2},\dots,\oray{v}{u_{-k}},\ray{v}{u_k})$.
\end{proof}

\Figure{BaseCase}{\includegraphics{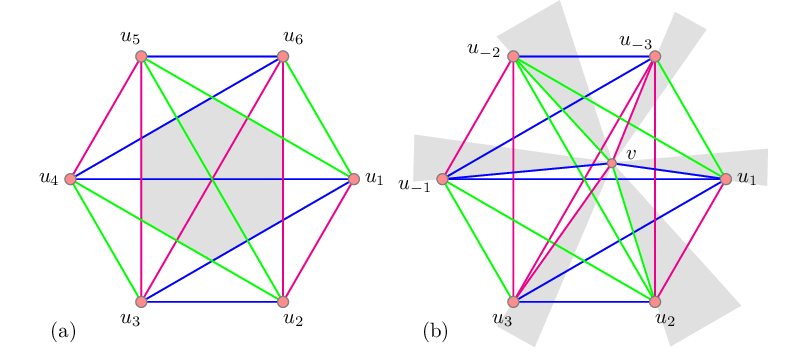}}{(a) Book embedding of $K_{2k}$. (b) Geometric drawing of $K_{2k+1}$ in which $v$ is a fan.}


The next proposition implies that $\gt{\T_{2k}}\leq k$, thus completing the proof of \thmref{GeometricThickness}.

\begin{proposition} 
\proplabel{GeometricThickness}
For all $k\geq2$, every $2k$-tree $G$ has a good drawing.
\end{proposition}

\begin{proof}
In this proof we repeatedly use two indices, \ii\  and \rr, whose ranges remain unchanged; in particular, $\ii\in I$ and $\rr\in\{1,2,\dots,k\}$.

We proceed by induction on $|V(G)|$. If $G\simeq K_{2k+1}$ the result is \lemref{BaseCase}. Now assume that $G\not\simeq K_{2k+1}$. Apply \lemref{Partition} to $G$. We obtain a nonempty independent set $S$ of $2k$-simplicial vertices of $G$. 

First suppose that $G\setminus S\simeq K_{2k}$. Let $v$ be an arbitrary vertex in $S$. By \lemref{BaseCase}, $G\setminus(S\setminus\{v\})$ ($\simeq K_{2k+1}$) has a good drawing in which $v$ is a fan. By \lemref{Insert}, each vertex $w\in S\setminus\{v\}$ can be inserted into the drawing (one at the time) resulting in a good drawing of $G$. 

Otherwise, by \lemref{Partition}(b), $G\setminus S$ is a $2k$-tree containing a $2k$-simplicial vertex $v$, such that $N_G(w)\subset N_{G\setminus S}[v]$ for each vertex $w\in S$. 

Apply the induction hypothesis to $G\setminus S$. If $G\setminus S\simeq K_{2k+1}$ then we can nominate $v$ to be a vertex of $G\setminus S$ that is a fan. By induction, we obtain a good drawing $\phi$ of $G\setminus S$ in which $v$ is a fan. Say $N_{G\setminus S}(v)= (u_1,u_2,\dots,u_k,u_{-1},u_{-2},\dots,u_{-k})$ in clockwise order about $v$. Thus each edge $vu_\ii$ is coloured $|\ii|$. 

By \lemref{Partition}(b), for each vertex $w\in S$, there is exactly one $\ii\in I$ for which $N_G(w)=N_{G\setminus S}[v]\setminus\{u_\ii\}$. Let $S_\ii:=\{w\in S:N_G(w)=N_{G\setminus S}[v]\setminus\{u_\ii\}\}$ for each $\ii\in I$. The vertices in $S_\ii$ have the same neighbourhood in $G$, and $\{S_\ii:\ii\in I\}$ is a partition of $S$. 

For each $\ii\in I$,  choose one vertex $x_\ii\in S_\ii$ (if any). Let $Q:=\{x_\ii:\ii\in I\}$. Suppose we have a good drawing of $(G\setminus S)\cup Q$. Then  by \lemref{Insert}, each vertex $w\in S\setminus Q$ can be inserted into the drawing (one at the time) resulting in a good drawing of $G$. Thus, from now on, we can assume that $S=Q\ (=\{x_\ii:\ii\in I\})$. Below we describe how to insert the vertices $\{x_\ii:\ii\in I\}$ into $\phi$ to obtain a good drawing $\phi'$ of $G$.

First we colour the edges incident to each vertex $x_\ii\in S$. Colour $x_\ii v$ by $|\ii|$, and colour $x_\ii u_\jj$ by $|\jj|$ for all $\jj\in I\setminus\{i\}$. Thus there are exactly two edges of each colour incident to $x_\ii$. In particular, $x_\ii v$ and $x_\ii u_{-\ii}$ are coloured $|\ii|$, and $x_\ii u_\jj$ and $x_\ii u_{-\jj}$ are coloured $\jj$ for each $\jj\in\{1,2,\dots,k\}\setminus\{|\ii|\}$.

For each $\rr \in\{1,2\dots,k\}$, let $G_\rr$ be the spanning subgraph of $G$ consisting of all the edges of $G$ coloured $\rr$. Let $G^{\star}$ be the spanning subgraph of $G$ with edge set $E(G)\setminus\{vu_\ii:\ii\in I\}$. Let $G^{\star}_\rr:=G_\rr\cap G^{\star}$ for each $\rr\in\{1,2,\dots,k\}$. 

As an intermediate step, we now construct a geometric drawing $\phi^\star$ of $G^\star$ (not in general position), in which each subgraph $G_\rr^\star$ is noncrossing. We later modify $\phi^\star$, by moving each vertex $x_\ii$ and drawing each edge $vu_\ii$, to obtain a general position geometric drawing $\phi'$ of $G$, in which each subgraph $G_\rr$ is noncrossing.


First, let $\phi^\star(w):=\phi(w)$ for every vertex $w$ of $G\setminus S$. By \obsref{EmptyDisc}, $v$ is $\varepsilon$-empty in $\phi$ for some $\varepsilon>0$. We now position each vertex $x_\ii$ on the segment $\seg{v}{u_{-\ii}}\cap D_\varepsilon(v)$. We have $F_\jj(x_\ii)=\WEDGE{\ray{x_\ii}{u_\jj}}{\oray{x_\ii}{u_{-\jj}}}$ for all $\jj\in I\setminus\{\ii,-\ii\}$. Observe that with
$x_\ii\in\seg{v}{u_{-\ii}}\cap D_\varepsilon(v)$, we have $v\not\in F_\jj(x_\ii)$ for all $\jj\in I\setminus\{\ii, -\ii\}$. Therefore, for $\ii\in I$ in some arbitrary order, each vertex $x_\ii$ can be positioned on the segment $\seg{v}{u_{-\ii}}\cap D_\varepsilon(v)$ so that:\\
\POINT\ $x_\ii\not\in F_\jj(x_\ell)$ for each $\ell\in I\setminus \{\ii\}$ and $\jj\in I\setminus\{\ell, -\ell\}$, and \\
\POINT\ $V(G)$ is in general position except for the collinear triples $v,x_\ii,u_{-\ii}$ ($\ii\in I$).

This is possible by the previous observation, since there is always a point close enough to $v$ where $x_\ii$ can be positioned. This placement of vertices of $G^\star$ determines a geometric drawing $\phi^\star$ of $G^\star$. The construction is illustrated in \figref{NewConstruction}. 

\begin{claim}
The subgraph $G^\star_\rr$ is noncrossing in $\phi^\star$ for each $\rr\in\{1,2,\dots,k\}$.
\end{claim}

\begin{proof}
Distinguish the following three types of edges in $G^\star_\rr$: \\
\ONE\ edges of $G_\rr^\star\setminus S$, \\
\TWO\ the edges $x_\rr v$, $x_\rr u_{-\rr}$, $x_{-\rr}v$, and $x_{-\rr}u_\ii$,\\
\THREE\ edges $x_\jj u_\rr$ and $x_\ell u_{-\rr}$ for distinct $\jj,\ell\in I\setminus\{\rr, -\rr\}$.

First note that $\segc{x_\rr}{v}\cup \segc{x_\rr}{u_{-\rr}}=\segc{v}{u_{-\rr}}$, and similarly $\segc{x_{-\rr}}{v}\cup \segc{x_{-\rr}}{u_\rr}=\segc{v}{u_\rr}$. Since no pair of edges in $\{vu_{-\rr}, vu_\rr\}\cup E(G_\rr\setminus S)$ cross in $\phi$, no pair of edges in $\{x_\rr v, x_\rr u_{-\rr}, x_{-\rr}v, x_{-\rr}u_\rr\} \cup E(G_\rr\setminus S)$ cross in $\phi^\star$.

It remains to prove that no type-(1) edge crosses a type-(3) edge, no type-(2) edge crosses a type-(3) edge, and that no pair of type-(3) edges cross in $\phi^\star$.

\Figure[H]{NewConstruction}{\hspace*{-2mm}\vspace*{-1mm}\includegraphics{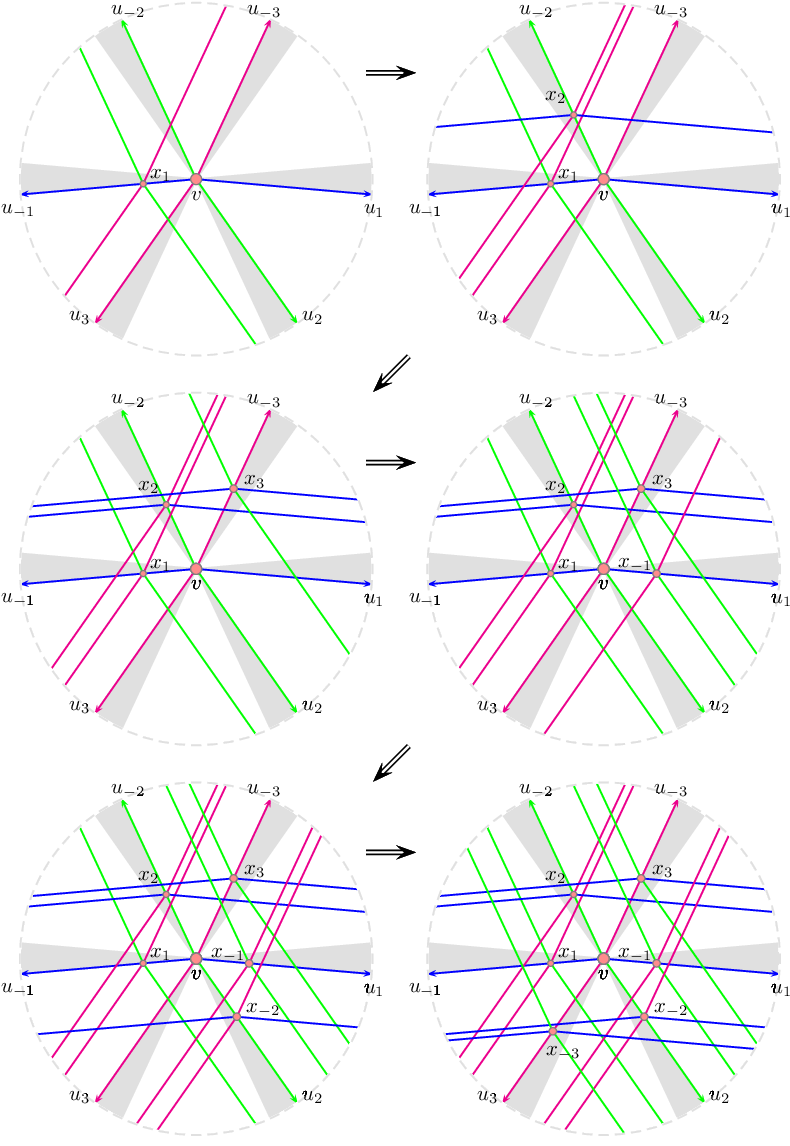}}{Placing each $x_\ii$ on the segment \seg{v}{u_{-\ii}}; intuitively speaking, the circle $D_\varepsilon$ is chosen small enough so that the edges incident with $u_\ii$ are almost parallel.}

Consider a type-(1) edge $e$ and a type-(3) edge. Since $v$ is a fan in $\phi$, the only two edges coloured \rr\ that are incident to $v$ in $G\setminus S$ (and thus in $G_\rr^\star\setminus S$) are $vu_\rr$ and $vu_{-\rr}$. These two edges are not in $G_\rr^\star$, and thus $e$ is not incident to $v$. Then, since $x_\jj\in D_\varepsilon(v)$ for all $\jj\in I$, \obsref{MustCross} implies that $e$ crosses $\segc{v}{u_\rr}$ or $\segc{v}{u_{-\rr}}$. That is ruled out in the previous case (when considering type-(1) and type-(2) edges) since $\segc{v}{u_\rr}=\segc{x_{-\rr}}{v}\cup\ \segc{x_{-\rr}}{u_{\rr}}$  and $\segc{v}{u_{-\rr}}=\segc{x_\rr}{v}\cup\ \segc{x_\rr}{u_{-\rr}}$.

Now suppose that a type-(2) edge crosses a type-(3) edge $x_\jj u_\rr$. The edge $x_\jj u_\rr$ shares an endpoint with the segment $\segc{v}{u_\rr}$; thus $x_\jj u_\rr$  crosses neither $x_{-\rr}v$ nor $x_{-\rr}u_\rr$. If $x_\jj u_\rr$ crosses $\segc{v}{u_{-\rr}}$, then $v\in F_\rr(x_\jj)$, contradicting our placement of $x_\jj$. Thus $x_\jj u_\rr$ crosses $\segc{v}{u_{-\rr}}$ and therefore crosses neither $x_{\rr}v$ nor $x_{\rr}u_{-\rr}$. By symmetry, no type-(2) edge crosses a type-(3) edge $x_{\ell}u_{-\rr}$.
 
Finally suppose that a type-(3) edge $x_\jj u_\rr$ crosses a type-(3) edge $x_\ell u_{-\rr}$. Then $x_\ell\in F_\rr(x_\jj)$ and $x_\jj\in F_{-\rr}(x_\ell)$, contradicting our placement of $x_\ell$ or $x_\jj$. Thus two type-(3) edges do not cross.
\end{proof} 

This completes the proof that $\phi^\star$ is a geometric drawing of $G^{\star}$, in which each $G^{\star}_\rr$ is noncrossing. The only collinear vertices in $\phi^\star$ are $v,x_\ii,u_{-\ii}$ for $\ii\in I$. We now move each vertex $x_\ii$ off the segment $\seg{v}{u_{-\ii}}\cap D_\varepsilon(v)$ into the wedge $F_{-\ii}(v)\cap D_\varepsilon(v)$. We achieve that with the help of \lemref{NewLemma} in the Appendix. 

For each $\rr\in\{1,2, \dots,k\}$, apply \lemref{NewLemma} to the noncrossing drawing of $G_\rr^{\star}$ induced by $\phi^{\star}$. We obtain a number $\varepsilon_\rr>0$ such that if $\phi':V(G_\rr^\star)\rightarrow\R^2$ is an injection with $\phi'(w)\in D_{\varepsilon_\rr}(\phi(w))$ for every vertex $w\in V(G^\star_\rr)$, then $\phi'$ is a noncrossing geometric drawing of $G^\star_\rr$ with the property that if three vertices $\phi'(a),\phi'(b),\phi'(c)$ are collinear in $\phi'$, then $\phi^\star(a),\phi^\star(b),\phi^\star(c)$ are collinear in $\phi^\star$. 

Let $\delta:=\min\{\varepsilon_\rr:\rr\in\{1,2, \dots,k\}\}$. For each $\ii\in I$, let $\phi'(x_\ii)$ be some point in the region $D_\delta(x_\ii)\cap F_{-\ii}(v)$. Let $\phi'(w):=\phi(w)$ for every other vertex of $G$. 

We now prove that each subgraph $G_\rr$ is noncrossing in $\phi'$. By \lemref{NewLemma}, since $\delta\leq\varepsilon_\rr$, each subgraph $G_\rr^{\star}$ is noncrossing in $\phi'$. We must also show that the edges $vu_\rr$ and $vu_{-\rr}$ do not cross any edge in $G_\rr$. First note that $vu_\rr$ and $vu_{-\rr}$ do not cross since they have a common endpoint. Suppose that an edge $e$ of $G_\rr^{\star}$ crosses $vu_{-\rr}$. Since the interior of the triangle $vx_\rr u_{-\rr}$ contains no vertex, $e$ also crosses $vx_\rr$ or $x_\rr u_{-\rr}$. This is impossible, since $vx_\rr$ and $x_\rr u_{-\rr}$ are edges of $G^{\star}_\rr$. Similarly, an edge $e$ of $G_\rr^{\star}$ does not cross $vu_\rr$. Thus $G_\rr$ is noncrossing. 

We now prove that $\phi'$ is in general position. By \lemref{NewLemma}, if three vertices are collinear in $\phi'$ then they are collinear in $\phi^\star$. The only collinear triples in $\phi^\star$ are $v,x_\ii,u_{-\ii}$ for $\ii\in I$. Since $\phi'(x_\ii)$ is in (the interior) of $F_{-\ii}(v)$, the vertices  $v,x_\ii,u_{-\ii}$ are not collinear in $\phi'$. Thus $\phi'$ is in general position. 

It remains to prove that every $2k$-simplicial vertex of $G$ is a fan in $\phi'$. Consider a $2k$-simplicial vertex $y$ that is not in $S$. By \lemref{Partition}(b), $y$ is not adjacent to $v$. Thus $y$ is adjacent to no vertex in $S$, and $y$ is $2k$-simplicial in $G\setminus S$. Moreover, $G\setminus S$ is not complete. By induction, $y$ is a fan in the drawing of $G\setminus S$ induced by $\phi'$, and thus $y$ is a fan in $\phi'$. Each vertex in $S$ is a fan in $\phi'$ by the following claim.

\begin{claim} 
For each $\ii\in I$, the vertex $x_\ii$ is a fan in $\phi'$.
\end{claim}

\begin{proof}
Let $H$ be the $2k$-tree obtained from $G\setminus S$ by adding a new vertex $h$ onto the $2k$-clique $N_{G\setminus S}(v)=\{u_\ii,u_{-\ii}:\ii\in I\}$. Consider the general position geometric drawing $\sigma$ of $H$ induced by $\phi$ with $\sigma(h):=\phi'(x_\ii)$. 


By construction, $\sigma(h)\in D_\varepsilon(v)\cap F_{-\ii}(v)$. Thus property (d) of the choice of $\varepsilon$ implies that the clockwise orders of $R(v,N_{G\setminus S}(v))$ and $R(h,N_H(h))$ are the same. Since $v$ is balanced in $\phi$, $h$ is balanced in $\sigma$. 

Now consider the drawings $\phi'$ of $G$ and $\sigma$ of $H$. Note that $F_\jj(x_\ii)=F_\jj(h)$ for all $\jj\in I\setminus\{\ii,-\ii\}$. Furthermore, since $\lin{x_\ii}{v}\subset F_{\ii}(h)\cup F_{-\ii}(h)\cup\{h\}$, we have $F_\ii(x_\ii)\subset F_\ii(h)$ and $F_{-\ii}(x_\ii)\subset F_{-\ii}(h)$. Therefore, $x_\ii$ is balanced in $\phi'$. 
The edges $x_\ii v$ and $x_\ii u_{-i}$ are both coloured $|\ii|$, and $x_\ii u_\jj$ is coloured $|\jj|$ for all $\jj\in I\setminus\{\ii\}$. Therefore $x_\ii$ is a fan in $\phi'$.
\end{proof}

We have thus proved that $\phi'$ is a general position geometric drawing of $G$, such that for each $\rr\in\{1,2,\dots, k\}$, the induced drawing of $G_\rr$ is noncrossing, and every $2k$-simplicial vertex is a fan. Thus $\phi'$ has thickness $k$ and is a good drawing of $G$. This completes the proof of \propref{GeometricThickness}, which implies \thmref{GeometricThickness}. 
\end{proof}

Note that it is easily seen that each noncrossing subgraph $G_\rr$ in the proof of \propref{GeometricThickness} is series-parallel.

\section{Book Thickness Lower Bound}
\seclabel{BookThickness}

Here we prove \thmref{BookThickness} for $k\geq3$. By the discussion in \secref{Results}, it suffices to construct a $k$-tree $Q_k$ with book thickness $\bt{Q_k}\geq k+1$ for all $k\geq 3$. To do so, start with the $k$-tree $K^\star_{k,2k^2+1}$ defined in \secref{AbstractResults}. Recall that $K$ is a $k$-clique and $S$ is a set of $2k^2+1$ $k$-simplicial vertices in  $K^\star_{k,2k^2+1}$. For each vertex $v\in S$, choose three distinct vertices $x_1,x_2,x_3\in K$, and for each $i\in\{1,2,3\}$, add a set of four vertices onto the $k$-clique $(K\cup\{v\})\setminus\{x_i\}$. Each set of four vertices is called an \emph{$i$-block of} $v$. Let $T$ be the set of vertices added in this step. Clearly $Q_k$ is a $k$-tree; see \figref{LowerBound}.

\Figure{LowerBound}{\includegraphics{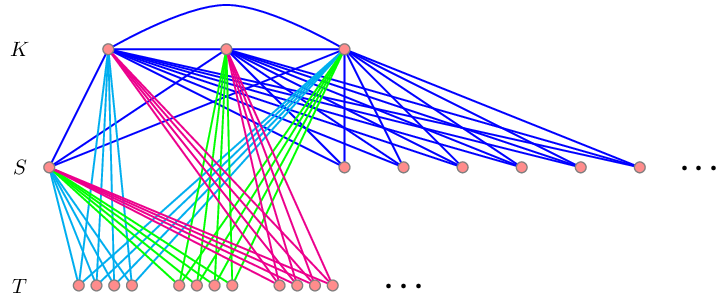}}{The graph $Q_k$ in \thmref{BookThickness} with $k=3$.}

\begin{lemma}
\lemlabel{LowerBound}
The book thickness of $Q_k$ satisfies $\bt{Q_k}\geq k+1$.
\end{lemma}

\begin{proof}
Suppose, for the sake of contradiction, that $Q_k$ has a book embedding with thickness $k$. Let $\{E_1,E_2,\dots,E_k\}$ be the corresponding partition of the edges. For each ordered pair of vertices $v,w\in V(Q_k)$, let the \emph{arc-set} $\as{vw}$ be the list of vertices in clockwise order from $v$ to $w$ (not including $v$ and $w$).

Say $K=(u_1,u_2,\dots,u_k)$ in anticlockwise order. There are $2k^2+1$ vertices in $S$. Without loss of generality there are at least $2k+1$ vertices in $S\cap\as{u_1u_k}$. Let $(v_1,v_2,\dots,v_{2k+1})$ be $2k+1$ vertices in $S\cap\as{u_1u_k}$ in clockwise order.

Observe that the $k$ edges $\{u_iv_{k-i+1}:1\leq i\leq k\}$ are pairwise crossing, and thus receive distinct colours, as illustrated in \figref{LowerBoundProof}(a). Without loss of generality, each $u_iv_{k-i+1}\in E_i$. As illustrated in \figref{LowerBoundProof}(b), this implies that $u_1v_{2k+1}\in E_1$, since $u_1v_{2k+1}$ crosses all of $\{u_iv_{k-i+1}:2\leq i\leq k\}$ which are coloured $\{2,3,\dots,k\}$. As illustrated in \figref{LowerBoundProof}(c), this in turn implies $u_2v_{2k}\in E_2$, and so on. By an easy induction, we obtain that $u_iv_{2k+2-i}\in E_i$ for all $i\in\{1,2,\dots,k\}$, as illustrated in \figref{LowerBoundProof}(d). It follows that for all $i\in\{1,2,\dots,k\}$ and $j\in\{k-i+1,k-i+2,\dots,2k+2-i\}$, the edge $u_iv_j\in E_i$, as illustrated in \figref{LowerBoundProof}(e). Finally,
as illustrated in \figref{LowerBoundProof}(f), we have:

\smallskip
If $qu_i\in E(Q_k)$ and $q\in\as{v_{k-1}v_{k+3}}$, then $qu_i\in E_i$.\hfill ($\star$)
\smallskip

\Figure{LowerBoundProof}{\includegraphics{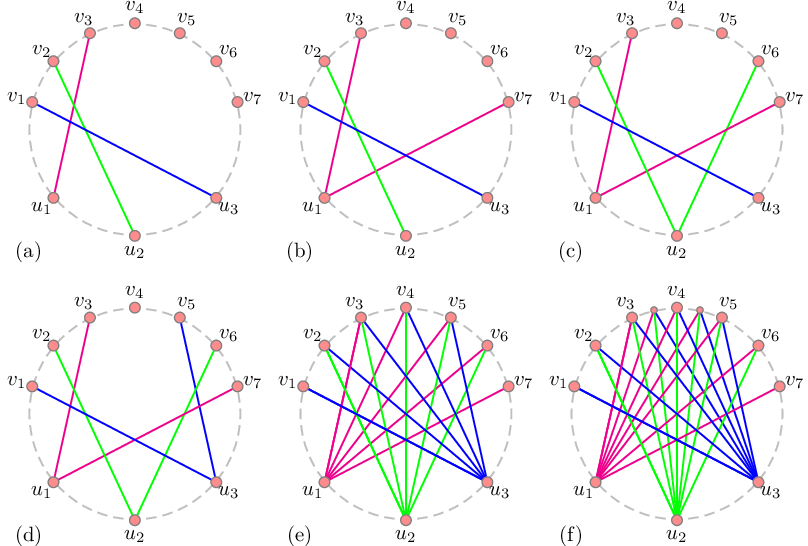}}{Illustration of the proof of \lemref{LowerBound} with $k=3$.}

Consider one of the twelve vertices $w\in T$ that are added onto a clique that contain $v_{k+1}$. Then $w$ is adjacent to $v_{k+1}$. Moreover, $w$ is in $\as{v_{k}v_{k+1}}$ or $\as{v_{k+1}v_{k+2}}$, as otherwise the edge $wv_{k+1}$ crosses $k$ edges of $Q_k[\{v_{k-1},v_{k+1}\};K]$ that are all coloured differently, which is a contradiction. By the pigeon-hole principle, one of $\as{v_{k}v_{k+1}}$ and $\as{v_{k+1}v_{k+2}}$ contains at least two vertices from two distinct $p$-blocks of $v_{k+1}$. Without loss of generality, $\as{v_kv_{k+1}}$ does. Let these four vertices be $(a,b,c,d)$ in clockwise order. 

Each vertex in $\{b,c,d\}$ is adjacent to $k-1$ vertices of $K$. Not all of $b,c,d$ are adjacent to the same subset of $k-1$ vertices in $K$, as otherwise all of $b,c,d$ would belong to the same $p$-block. Hence each vertex in $K$\ has a neighbour in $\{b,c,d\}$. By ($\star$) the edges of $Q_k[\{b,c,d\},K]$ receive all $k$ colours. However, every edge in $Q_k[\{b,c,d\};K]$ crosses the edge $av_{k+1}$, implying that there is no colour available for $av_{k+1}$. 
This contradiction completes the proof.
\end{proof}
 
Note that the number of vertices in $Q_k$ is $|K|+|S|+|T|=k+2k^2+1+3\cdot4\cdot(2k^2+1)=13(2k^2+1)+k$. Adding more simplicial vertices to $Q_k$ cannot reduce its book thickness. Thus for all $n\geq13(2k^2+1)+k$, there is a $k$-tree $G$ with $n$ vertices and $\bt{G}=k+1$.

\section{Open Problems}
\seclabel{Problems}

\paragraph{\bf Complete Graphs:} The thickness of the complete graph $K_n$ was intensely studied in the 1960's and 1970's. Results by a number of authors \citep{AG76, Beineke67, BH65, Mayer72} together prove that  $\tt{K_n}=\ceil{(n+2)/6}$, unless $n=9$ or $10$, in which case $\tt{K_9}=\tt{K_{10}}=3$. \citet{BK79} proved that $\bt{K_n}=\ceil{n/2}$. In fact, it is easily seen that
\begin{equation*}
\arb{K_n}=\ga{K_n}=\bt{K_n}=\ba{K_n}=\ceil{n/2}\enspace.
\end{equation*}
\citet{AK85} proved that $\sa{K_n}=\ceil{n/2}+1$. Now
\begin{equation*}
\gsa{K_n}\leq\bsa{K_n}\leq n-1\enspace.
\end{equation*}
(\emph{Proof}. Place the vertices of $K_n$ on a circle, with a spanning star rooted at each vertex except one.)\ What is \gsa{K_n} and \bsa{K_n}?

\citet{BHRW-CGTA06} proved that every geometric drawing of $K_n$ has arboricity (and thus thickness) at most $n-\sqrt{n/12}$. It is unknown whether for some constant $\varepsilon>0$, every geometric drawing of $K_n$ has thickness at most $(1-\varepsilon)n$; see \citep{BHRW-CGTA06}. \citet{DEH-JGAA00} studied the geometric thickness of $K_n$, and proved that\footnote{\citet{Archdeacon-GeomThickness} writes, ``The question (of the value of $\gt{K_n}$) was apparently first raised by Greenberg in some unpublished work. I read some of his personal notes in the library of the University of Riga in Latvia. He gave a construction that showed $\gt{K_n}\leq\ceil{n/4}$.''} 
\begin{equation*}
\ceil{(n/5.646)+0.342}\leq\gt{K_n}\leq\ceil{n/4}\enspace.
\end{equation*}
What is \gt{K_n}? It seems likely that the answer is closer to $\ceil{n/4}$ rather than the above lower bound.

\paragraph{\bf Asymptotics:} \citet{Eppstein01} (also see
\citep{Blankenship-PhD03}) constructed $n$-vertex graphs $G_n$ with
$\sa{G_n}=\arb{G_n}=\tt{G_n}=\gt{G_n}=2$ and $\bt{G_n}\rightarrow\infty$. Thus
book thickness is not bounded by any function of geometric thickness. Similarly,
\citet{Eppstein-AMS} constructed $n$-vertex graphs $H_n$ with
$\sa{H_n}=\arb{H_n}=\tt{H_n}=3$ and $\gt{H_n}\rightarrow\infty$. Thus geometric
thickness is not bounded by any function of thickness (or arboricity).
\citet{Eppstein-AMS} asked whether graphs with thickness $2$ have bounded
geometric thickness? Whether all graphs with arboricity $2$ have bounded geometric thickness is also interesting. It is easily seen that graphs with star arboricity $2$ have geometric star arboricity at most $2$ (cf.~\citep{Brass-etal-CGTA07}).

\paragraph{\bf Book Arboricity:} \citet{BK79} proved that every graph $G$ with book thickness $t$ satisfies $|E(G)|\leq(t+1)|V(G)|-3t$. Thus \Eqnref{Arboricity} implies that $\arb{G}\leq\bt{G}+1$ for every graph $G$, as observed by \citet{DH91}. Is $\ba{G}\leq\bt{G}+1$?

\paragraph{\bf Number of Edges:} Let $\E_m$ be the class of graphs with at most $m$ edges. \citet{DHS91} proved that $\tt{\E_m}\leq\sqrt{m/3}+3/2$. What is the minimum $c$ such that $\tt{\E_m}\leq(c+o(1))\sqrt{m}$? \citet{DHS91} conjectured that the answer is $c=1/16$, which would be tight for the balanced complete bipartite graph \citep{BHM64}. \citet{Malitz94a} proved using a probabilistic argument that $\bt{\E_m}\leq72\sqrt{m}$. Is there a constructive proof that $\bt{\E_m}\in\Oh{\sqrt{m}}$ or $\gt{\E_m}\in\Oh{\sqrt{m}}$? What is the minimum $c$ such that $\gt{\E_m}\leq(c+o(1))\sqrt{m}$ or $\bt{\E_m}\leq(c+o(1))\sqrt{m}$? 

\paragraph{\bf Planar Graphs:} Recall that \citet{Yannakakis89} proved that every planar graph $G$ has book thickness $\bt{G}\leq 4$. He also claims there is a planar graph $G$ with $\bt{G}=4$. A construction is given in the conference version of his paper \citep{Yannakakis86b}, but the proof is far from complete: Yannakakis admits, ``Of course, there are many other ways to lay out the graph'' \citep{Yannakakis86b}. The journal version \citep{Yannakakis89} cites a paper ``in preparation'' that proves the lower bound. This paper has not been published. Therefore we consider it an open problem whether $\bt{G}\leq3$ for every planar graph $G$. 

Let $G_0=K_3$. For $k\geq1$, let $G_k$ be the planar $3$-tree obtained by adding a $3$-simplicial vertex onto the vertex set of each face of $G_{k-1}$. In the journal version of this paper we conjectured that $\bt{G_k}=4$ for sufficiently large $k$. This conjecture is false. Indeed, Heath [Embedding planar graphs in seven pages, \emph{Proc. 25th Annual Symp. on Foundations of Comput. Sci.} (FOCS '84), 74--83, IEEE, 1984] proved that every graph with treewidth 3 (which includes $G_k$) has a 3-page book embedding. \label{CHANGE} Note that Heath's proof implicitly uses the notion of colourful simplicial vertices (which we reinvented). 

\paragraph{\bf Genus:} Let $\X_\gamma$ denote the class of graphs with genus at most $\gamma$. \citet{DH91} proved that $\tt{\X_\gamma}\leq 6+\sqrt{2\gamma-2}$; also see \citep{Asano-AC94, Asano-JCTB87}. What is the minimum $c$ such that $\tt{\X_\gamma}\leq(c+o(1))\sqrt{\gamma}$? Building on prior work by \citet{HI92}, \citet{Malitz94b} proved using a probabilistic argument that $\bt{\X_\gamma}\in\Oh{\sqrt{\gamma}}$, and thus $\gt{\X_\gamma}\in\Oh{\sqrt{\gamma}}$. Is there a constructive proof that $\bt{\X_\gamma}\in\Oh{\sqrt{\gamma}}$ or $\gt{\X_\gamma}\in\Oh{\sqrt{\gamma}}$. What is the minimum $c$ such that $\bt{\X_\gamma}\leq(c+o(1))\sqrt{\gamma}$, or $\gt{\X_\gamma}\leq(c+o(1))\sqrt{\gamma}$? 

\citet{Endo97} proved that $\bt{\X_1}\leq7$. Let $\chi(\X_\gamma)$ denote the maximum chromatic number of all graphs with genus at most $\gamma$. Heawood's formula and the four-colour theorem state that $\chi(\X_\gamma)=\floor{\half(7+\sqrt{1+48\gamma})}$. Thus $\chi(\X_\gamma)$ and the known upper bounds on $\bt{\X_\gamma}$ coincide for $\gamma=0$ and $\gamma=1$. \citet{Endo97} asked whether $\bt{\X_\gamma}=\chi(\X_\gamma)$ for all $\gamma$. Both $\bt{\X_\gamma}$ and $\chi(\X_\gamma)$ are in $\Oh{\sqrt{\gamma}}$. There is some tangible evidence relating book thickness and chromatic number. First, \citet{BK79} proved that $\chi(G)\leq2\cdot\bt{G}+2$ for every graph $G$. Second, the maximum book thickness and maximum chromatic number coincide ($=k+1$) for graphs of treewidth $k\geq3$. In fact, the proof by \citet{GH-DAM01} that $\bt{\T_k}\leq k+1$ is based on the ($k+1$)-colourability of $k$-trees.

\paragraph{\bf Minors:} Let $\M_\ell$ be the class of graphs with no $K_\ell$-minor. Note that $\M_3=\T_1$ and $\M_4=\T_2$. \citet{JMOS98} proved that $\tt{\M_5}=2$. What is $\gt{\M_5}$ and $\bt{M_5}$? \citet{Kostochka82} and \citet{Thomason84} independently proved that the maximum arboricity of all graphs with no $K_\ell$ minor is $\Theta(\ell\sqrt{\log\ell})$. In fact, \citet{Thomason01} asymptotically determined the right constant. Thus $\tt{\M_\ell}\in\Theta(\ell\sqrt{\log\ell})$ by \Eqnref{AbstractRelationship}. Blankenship and Oporowski \citep{Blankenship-PhD03,BO01} proved that $\bt{\M_\ell}$ (and hence $\gt{\M_\ell}$) is finite. The proof depends on Robertson and Seymour's deep structural characterisation of the graphs in $\M_\ell$. As a result, the bound on $\bt{\M_\ell}$ is a truly huge function of $\ell$. Is there a simple proof that $\gt{\M_\ell}$ or $\bt{\M_\ell}$ is finite? What is the right order of magnitude of $\gt{\M_\ell}$ and $\bt{\M_\ell}$?

\paragraph{\bf Maximum Degree:} Let $\D_\Delta$ be the class of graphs with maximum degree at most $\Delta$. \citet{Wessel83} and \citet{Halton91} independently proved that $\tt{\D_\Delta}\leq\ceil{\Delta/2}$, and \citet{SSV-IS04} proved that $\tt{\D_\Delta}\geq\ceil{\Delta/2}$. Thus $\tt{\D_\Delta}=\ceil{\Delta/2}$. \citet{Eppstein-AMS} asked whether $\gt{\D_\Delta}$ is finite. A positive result in this direction was obtained by \citet{DEK-SoCG04}, who proved that $\gt{\D_4}\leq2$. On the other hand, \citet{BMW-EJC06} recently proved that $\gt{\D_\Delta}=\infty$ for all $\Delta\geq 9$; in particular, there exists $\Delta$-regular $n$-vertex graphs with geometric thickness $\Omega(\sqrt{\Delta}n^{1/2-4/\Delta-\varepsilon})$. It is unknown whether $\gt{\D_\Delta}$ is finite for $\Delta\in\{5,6,7,8\}$. 

\citet{Malitz94a} proved that there exists $\Delta$-regular $n$-vertex graphs with book thickness \linebreak
$\Omega(\sqrt{\Delta}n^{1/2-1/\Delta})$. \citet{BMW-EJC06} reached the same conclusion for all $\Delta\geq3$. Thus $\bt{\D_\Delta}=\infty$ unless $\Delta\leq2$. Open problems remain for specific values of $\Delta$. For example, the best bounds on $\bt{\D_3}$ are $\Omega(n^{1/6})$ and \Oh{n^{1/2}}.

\paragraph{\bf Computational Complexity:} Arboricity can be computed in polynomial time using the matroid partitioning algorithm of \citet{Edmonds65}. Computing the thickness of a graph is \NP-hard \citep{Mansfield83}. Testing whether a graph has book thickness at most $2$ is \NP-complete \citep{Wigderson82}. \citet{DEH-JGAA00} asked what is the complexity of determining the geometric thickness of a given graph? The same question can be asked for all of the other parameters discussed in this paper.

\subsection*{Acknowledgement} 

Thanks to the anonymous reviewer for numerous helpful comments.


\def\soft#1{\leavevmode\setbox0=\hbox{h}\dimen7=\ht0\advance \dimen7
  by-1ex\relax\if t#1\relax\rlap{\raise.6\dimen7
  \hbox{\kern.3ex\char'47}}#1\relax\else\if T#1\relax
  \rlap{\raise.5\dimen7\hbox{\kern1.3ex\char'47}}#1\relax \else\if
  d#1\relax\rlap{\raise.5\dimen7\hbox{\kern.9ex \char'47}}#1\relax\else\if
  D#1\relax\rlap{\raise.5\dimen7 \hbox{\kern1.4ex\char'47}}#1\relax\else\if
  l#1\relax \rlap{\raise.5\dimen7\hbox{\kern.4ex\char'47}}#1\relax \else\if
  L#1\relax\rlap{\raise.5\dimen7\hbox{\kern.7ex
  \char'47}}#1\relax\else\message{accent \string\soft \space #1 not
  defined!}#1\relax\fi\fi\fi\fi\fi\fi} \def\cprime{$'$}


\section*{Appendix}

Here we prove two perturbation lemmas that are used in the proof of \thmref{GeometricThickness}. Three discs are \emph{collinear} if there is a line that intersects each disc.

\begin{lemma}
\lemlabel{BlowUpDiscs}
Let $P$ be a finite set of points in the plane. Then there exists $\varepsilon>0$ such that\\
\aaa\ $D_\varepsilon(u)\cap D_\varepsilon(v)=\emptyset$ for all $u,v\in P$,\\
\bbb\ for all $u,v,w \in P$, if the discs $D_\varepsilon(u), D_\varepsilon(v), D_\varepsilon(w)$ are collinear, \\
\hspace*{1.8em} then the points $u,v,w$ are collinear.
\end{lemma}

\begin{proof}
Say $P=\{p_1,\dots,p_n\}$. We prove the following statement by induction on $\ell\in\{0,1,\dots,n\}$:

\medskip
\begin{minipage}{\textwidth-5mm} 
For all $i\in\{1,2,\dots,\ell\}$, there exists a disk $D^\ell(p_i)$ of positive radius centered at $p_i$ such that the following two properties hold, where $D^{\ell}(p_i):=\{p_i\}$ for each $i>\ell$:\\
\aaa\ $D^\ell(p_i)\cap D^\ell(p_j)=\emptyset$ for all $p_i, p_j\in P$,\\
\bbb\ for all $p_i, p_j, p_k\in P$, if the discs $D^\ell (p_i), D^\ell (p_j), D^\ell (p_k)$ are collinear,\\
\hspace*{1.8em} then the points $p_i$, $p_j$, $p_k$ are collinear.
\end{minipage}
\smallskip

\noindent This statement implies the lemma, by defining $\varepsilon$ to be the radius of the smallest disk $D^n(p_i)$.

The base case $\ell=0$ is vacuous. Now assume that $\ell>0$. For all distinct $i$ and $j$ such that $p_i,p_j,p_\ell$ are not collinear, every line that intersects $D^{\ell-1}(p_i)$ and $D^{\ell-1}(p_j)$ does not intersect $p_\ell$, by induction.  Thus $p_\ell$ is in an open region $R$ of the plane defined by the complement of the union of all such lines. Thus, there is an open disk $D\subset R$ of positive radius centered at $p_\ell$, such that $D$ does not intersect $D^{\ell-1}(p_i)$, for all $i\neq\ell$. Defining $D^\ell(p_\ell):=D$ and $D^\ell(p_i):=D^{\ell-1}(p_i)$ for all points $p_i\neq p_\ell$ completes the proof.
\end{proof}

\begin{lemma}
\lemlabel{NewLemma}
Let $\phi$ be a noncrossing geometric drawing of a graph $G$ \paran{not necessarily in general position}. Then there exists $\varepsilon>0$ such that if $\phi':V(G)\rightarrow\R^2$ is an injection with $\phi'(v)\in D_\varepsilon(\phi(v))$ for every vertex $v\in V(G)$, then $\phi'$ is a noncrossing geometric drawing of $G$ with the property that if three vertices $\phi'(u),\phi'(v),\phi'(w)$ are collinear in $\phi'$, then $\phi(u),\phi(v),\phi(w)$ are collinear in $\phi$.
\end{lemma}

\begin{proof}
Let $\varepsilon>0$ be the constant obtained by applying \lemref{BlowUpDiscs} to the point set $\{\phi(v):v\in V(G)\}$. We now prove that every function $\phi'$, as defined in the statement of the lemma, has the desired properties. 

Suppose there exists a line intersecting three vertices $\phi'(u),\phi'(v),\phi'(w)$ in $\phi'$. Then the same line passes through the discs $D_\varepsilon(\phi(u)),D_\varepsilon(\phi(v)), D_\varepsilon(\phi(w))$. Thus, by \lemref{BlowUpDiscs}, $\phi(u),\phi(v),\phi(w)$ are collinear in $\phi$.

It remains to prove that $\phi'$ is a noncrossing geometric drawing of $G$. For each vertex $v\in V(G)$, let $J(v):=D_\varepsilon(\phi(v))$.  Thus the image $\phi'(v)$ is in $J(v)$. For each edge $vw\in E(G)$, let $J(vw)$ be the region consisting of the union of all segments with one endpoint in $J(v)$ and the other endpoint in $J(w)$. Since $\phi'(v)\in J(v)$ and $\phi'(w)\in J(w)$, the image $\seg{\phi'(v)}{\phi'(w)}$ of the edge $vw$ is contained in $J(vw)$.

Thus to prove that $\phi'$ is a noncrossing geometric drawing of $G$, 
it suffices to prove that:\\
\qi\ $J(v)\cap J(xy)=\emptyset$, for every vertex $v\in V(G)$ and edge $xy\in E(G)$\\
\hspace*{1.8em} not incident to $v$, and \\
\qii\ $J(vw)\cap J(xy)=\emptyset$, for all edges $vw,xy\in E(G)$ with no\\
\hspace*{1.8em} common endpoint.

We now prove (i). First, suppose that $\phi(v),\phi(x),\phi(y)$ are 
collinear. Then without loss of generality, $\phi(x)$ is between 
$\phi(v)$ and $\phi(y)$, as otherwise $v$ intersects the edge $xy$ in 
$\phi$. Since $J(v)\cap J(x)=\emptyset$ by \lemref{BlowUpDiscs}, we have 
$J(v)\cap J(xy)=\emptyset$. Now suppose that $\phi(v),\phi(x),\phi(y)$ 
are not collinear. By \lemref{BlowUpDiscs}, $J(v),J(x),J(y)$ are not 
collinear, in which case $J(v)\cap J(xy)=\emptyset$.

We now prove (ii). Suppose on the contrary that $J(vw)\cap J(xy)\ne\emptyset$. First suppose that no three points in $\{\phi(v),\phi(w),\phi(x),\phi(y)\}$ are collinear. Then at least one of $vw$ and $xy$, say $vw$, lies on the convex hull of $\{\phi(v),\phi(w),\phi(x),\phi(y)\}$. Then $\lin{\phi(v)}{\phi(w)}$ does not intersect \seg{\phi(x)}{\phi(y)}. Therefore, the only way for $J(vw)\cap J(xy)\ne\emptyset$ is if $J(x)\cup J(v)\cap J(vw)\ne\emptyset$ or  $J(v)\cup J(w)\cap J(xy)\ne\emptyset$. This is impossible by \lemref{BlowUpDiscs}, since no three points in $\{\phi(v),\phi(w),\phi(x),\phi(y)\}$  are collinear. 

Now suppose that exactly three vertices in $\{\phi(v),\phi(w),\phi(x),\phi(y)\}$ are collinear. Without loss of generality, $\phi(v),\phi(w),\phi(x)$ are collinear and $\phi(w)$ is between $\phi(v)$ and $\phi(x)$. Then the only way for $J(vw)\cap J(xy)\ne\emptyset$ is if $J(w)\cap J(xy)\ne\emptyset$, which is 
impossible by \lemref{BlowUpDiscs}, since $\phi(w),\phi(x),\phi(y)$ are not collinear by assumption.

Finally assume that all four points $\phi(v),\phi(w),\phi(x),\phi(y)$  are collinear. Since $vw$ and $xy$ do not cross in $\phi$, we may assume without loss of generality, that $v,w,x,y$ are in this order on the line. Then $J(vw)\cap J(xy)\ne\emptyset$ only if  $J(w)\cap J(x)\ne\emptyset$, which is impossible by \lemref{BlowUpDiscs}.
\end{proof}
\end{document}